\documentclass[aos,MSNbibl,nameyear,dvips]{arximspdf}

\doi{10.1214/12-AOS1048} %kopijuoti is PTS
\volume{40}
\issue{6}
\pubyear{2012}
\firstpage{3050}
\lastpage{3076}

\makeatletter
\renewcommand{\citep}[1]{(\citeauthor{#1} \citeyear{#1})}
\newcommand{\eqref}[1]{(\ref{#1})}
\newtheorem{thmm}{Theorem}
\newtheorem{lem}{Lemma}
\newtheorem{cor}{Corollary}

\newcommand{\I}{\mathrm{I}}

\newcommand{\fatdot}{ \cdot}

\newcommand{\nats}{\mathbb{N}}
\newcommand{\real}{\mathbb{R}}

\newcommand{\opand}{\mbox{ and }}
\makeatother

\begin{document}
\begin{frontmatter}

\title{Variable transformation to obtain geometric ergodicity
in the random-walk Metropolis algorithm}
\runtitle{Morphisms for geometric ergodicity}

\begin{aug}
\author[A]{\fnms{Leif T.} \snm{Johnson}\corref{}\ead[label=e1]{ltjohnson@google.com}}
\and
\author[B]{\fnms{Charles J.} \snm{Geyer}\ead[label=e2]{charlie@stat.umn.edu}}
\runauthor{L.~T. Johnson and C.~J. Geyer}
\affiliation{Google Inc. and University of Minnesota}
\address[A]{Google Inc.\\
1600 Amphitheatre Parkway\\
Mountain View, California 94043\\
USA\\
\printead{e1}}

\address[B]{School of Statistics\\
University of Minnesota\\
313 Ford Hall \\
224 Church St. SE \\
Minneapolis, Minnesota 55455\\
USA\\
\printead{e2}}
\end{aug}

% HISTORY:
\received{\smonth{1} \syear{2012}}
\revised{\smonth{6} \syear{2012}}

% ABSTRACT
%
\begin{abstract}
A random-walk Metropolis sampler is geometrically ergodic if its
equilibrium density is super-exponentially light and satisfies
a curvature condition [\textit{Stochastic Process. Appl.} \textbf{85}
(2000) 341--361]. Many applications, including
Bayesian analysis with conjugate priors of logistic and Poisson regression
and of log-linear models for categorical data result in posterior distributions
that are not super-exponentially light. We show how to apply
the change-of-variable formula for diffeomorphisms to obtain new
densities that
do satisfy the conditions for geometric ergodicity. Sampling the new variable
and mapping the results back to the old gives a geometrically ergodic sampler
for the original variable. This method of obtaining geometric ergodicity
has very wide applicability.
\end{abstract}

% KEYWORDS
\begin{keyword}[class=AMS]
\kwd[Primary ]{60J05}
\kwd{65C05}
\kwd[; secondary ]{60J22}
\end{keyword}

\begin{keyword}
\kwd{Markov chain Monte Carlo}
\kwd{change of variable}
\kwd{exponential family}
\kwd{conjugate prior}
\kwd{Markov chain isomorphism}
\kwd{drift condition}
\kwd{Metropolis--Hastings--Green algorithm}
\end{keyword}

\end{frontmatter}

%s1 #&#
\section{Introduction}\label{sec1}

%%% Note: OUTSTANDING introductory sentence
Markov chain Monte Carlo (MCMC) using the Metro\-polis--Hastings--Green algorithm
[\citet{metropolis}, \citet{hastings}, \citet{green}] or its special
case the Gibbs sampler
[\citet{geman-geman}, \citet{tanner-wong}, \citet{gelfand-smith}]
has become very widely used
[\citet{gilks}, \citet{brooks}], especially after \citet
{gelfand-smith} pointed
out that
most Bayesian inference can be done using MCMC, and little can be done without
it.

In ordinary, independent and identically distributed Monte Carlo (OMC), the
asymptotic variance of estimates is easily calculated
[\citet{intro}, Section~1.7]. In MCMC, the properties of estimates are more
difficult to handle theoretically [\citet{intro}, Section~1.8]. A Markov
chain central limit theorem (CLT) may or may not hold
[\citet{tierney}, \citet{chan-geyer}]. If it does hold, the
asymptotic variance of
MCMC estimates is more difficult to estimate than for OMC estimates, but
estimating the asymptotic variance of the MCMC estimates is doable
[\citet{geyer-orange}, \citet{flegal-jones},
\citet{intro}, Section~1.10]. The CLT holds for all $L^{2 +
\varepsilon}$ functionals of a Markov chain if the Markov chain is
geometrically ergodic [\citet{chan-geyer}]. For a reversible Markov chain
[\citet{intro}, Section~1.5] the CLT holds for all $L^2$ functionals
if and
only if the Markov chain is geometrically ergodic
[\citet{roberts-rosenthal}]. The CLT may hold for some functionals of a
Markov chain when the Markov chain is not geometrically ergodic
[\citet{gordin-lifsic},
\citet{maigret},
\citet{kipnis-varadhan},
\citet{chan},
\citet{tierney},
\citet{chan-geyer},
\citeauthor{roberts-rosenthal} (\citeyear{roberts-rosenthal,roberts-rosenthal-survey}),
\citet{jones}],
but then it is usually very difficult to verify that a CLT exists for a
given functional of the Markov chain. Thus geometric ergodicity is a very
desirable property for a Markov chain to have. This is especially true
because most instances of the Metropolis--Hastings--Green algorithm are
reversible or can be made to be reversible [\citet{intro},
Sections~1.5, 1.12 and~1.17],
so, as stated above, geometric ergodicity implies the CLT
holds for all $L^2$ functionals of the Markov chain, which makes reversible
geometrically ergodic MCMC just as good as OMC in this respect.

Geometric ergodicity also plays a key role in the theory
of calculable nonasymptotic bounds for Markov chain estimators
[\citet{rosenthal},
\citet{latuszynski-niemiro},
{\L}atuszy{\'{n}}ski, Miasojedow and
Niemiro\break (\citeyear{latuszynski-miasojedow-niemiro})],
but is only half of what must be done to establish this type of result.
The other half is establishing a minorization condition.
The proof techniques involved in establishing geometric ergodicity
and in establishing minorization conditions, however, have little in common.
We deal only with establishing geometric ergodicity.

%s1.1 #&#
\subsection{The random-walk Metropolis algorithm}

The Metropolis--Hastings--Green algorithm generates a Markov chain
having a specified invariant probability distribution.
We restrict our attention to distributions of continuous random vectors,
those having a density $\pi$ with respect to Lebesgue measure on
$\real^k$.
If $\pi$ is only known up to a normalizing constant,
then the Metropolis--Hastings--Green algorithm still works.

We describe only the random-walk Metropolis algorithm
[terminology introduced by \citet{tierney}]. This
simulates a Markov chain $X_1, X_2, \ldots$ having $\pi$ as an invariant
distribution. It is determined by $\pi$
and another function $q \dvtx\real^k \to\real$ that is a properly normalized
probability density with respect to Lebesgue measure on $\real^k$ and
is symmetric about zero.
Each iteration does the following three steps, where
$X_n$ is the state of the Markov chain before the iteration
and $X_{n + 1}$ is the state after the iteration.
Simulate $Z_n$ having the distribution $q$, and set $Y_n = X_n + Z_n$.
Calculate
%
%
%e1 #&#
\begin{equation}
\label{eqaccept} a(X_n, Y_n) = \min \bigl(1,
\pi(Y_n) / \pi(X_n) \bigr).
\end{equation}
Set $X_{n+1} = Y_n$ with probability $a(X_n, Y_n)$,
and set $X_{n+1} = X_n$ with probability $1 - a(X_n, Y_n)$.

The only requirement is $\pi(X_1) > 0$. The operation of the algorithm
itself then ensures that $\pi(X_n) > 0$ almost surely for all $n$,
so \eqref{eqaccept} always makes sense.

The proposal density $q$ and target density $\pi$ are arbitrary.
The algorithm always produces a (not necessarily ergodic) reversible
Markov chain having invariant density $\pi$ regardless of what $q$ is chosen.
If $q$ is everywhere positive, then the Markov chain is necessarily ergodic
[irreducible and positive Harris recurrent, \citet{tierney}, Corollary~2].

The R package \texttt{mcmc} [\citet{mcmc-R-package}]
provides a user-friendly implementation of
the random-walk Metropolis algorithm combined with the variable transformation
methodology described in this article in
its \texttt{morph.metrop} function.
The user provides an R function that
evaluates $\log\pi$, and the \texttt{metrop} function in that package
does the simulation. If the user correctly codes the function
that evaluates $\log\pi$, then the \texttt{morph.metrop} function is
guaranteed
to simulate a reversible ergodic Markov chain having invariant density
$\pi$.
This gives an algorithm having an enormous range of application,
which includes all Bayesian inference for models with continuous parameters
and continuous prior distributions.
No other computer package known to us combines this range of application
with the correctness guarantees of the \texttt{mcmc} package,
which are as strong as can be made about arbitrary user-specified target
distributions.

%s1.2 #&#
\subsection{Geometric ergodicity and random-walk Metropolis}

A random-walk Metropolis sampler is not necessarily geometrically ergodic,
but its geometric ergodicity has received more attention
[\citet{mengersen-tweedie},
\citet{roberts-tweedie},
\citet{jarner-hansen}] than any other
MCMC sampler, except perhaps independence Metropolis--Hastings samplers,
also terminology introduced by \citet{tierney},
which are also studied in \citet{mengersen-tweedie} and \citet
{roberts-tweedie}.
Independence Metropolis--Hastings samplers, however, do not have good
properties,
being either uniformly ergodic or not geometrically ergodic and uniformly
ergodic only when its proposal distribution is particularly adapted to
$\pi$
in a way that is difficult to achieve (whenever independence samplers
work, importance sampling also works, so MCMC is unnecessary).

To simplify the theory, \citet{mengersen-tweedie},
\citet{roberts-tweedie} and \citet{jarner-hansen} restrict attention to
$\pi$ that are strictly positive and continuously differentiable.
In order to build on their results, we also adopt this
restriction.\vadjust{\goodbreak}
The geometric ergodicity properties of the random-walk Metropolis
algorithm are related to
%
%
%e2 #&#
\begin{equation}
\label{eqtails} \limsup_{\vert x \vert \rightarrow\infty} \frac
{x}{\vert x \vert} \cdot\nabla\log\pi(x),
\end{equation}
where the dot indicates inner product, and $\vert\fatdot \vert$
denotes the
Euclidean norm.
We say~$\pi$ is \textit{super-exponentially light} if \eqref{eqtails}
is $- \infty$, is \emph{exponentially light} if \eqref{eqtails} is negative
and \emph{sub-exponentially light} if \eqref{eqtails} is zero.

None of these conditions are necessary for geometric ergodicity.
A necessary condition for the geometric ergodicity of a random-walk Metropolis
algorithm is that the target density $\pi$ have a moment generating function
[\citet{jarner-tweedie}].
It is possible for a density to have a moment generating function but
not be
even sub-exponentially light, for example, the unnormalized density
\[
\pi(x) = e^{- \vert x \vert} \bigl(1 + \cos(x) \bigr),\qquad x \in \real.
\]

Following \citet{roberts-tweedie} and \citet{jarner-hansen}, we also
restrict attention to $q$ that are bounded away from zero in a neighborhood
of zero. This includes the normal proposal distributions used by the
R package \texttt{mcmc}.
%
%
%th1 #&#
\begin{thmm}[{[\protect{\citet{jarner-hansen}, Theorem 4.3}]}]\label{thmjh-main}
Suppose $\pi$ is a super-exponentially light density on $\real^k$ that
also satisfies
%
%
%e3 #&#
\begin{equation}
\label{condjh-curvature} \limsup_{\vert x \vert \rightarrow\infty} \frac{x}{\vert x \vert} \cdot
\frac{\nabla\pi(x)}{\vert\nabla\pi(x) \vert} < 0,
\end{equation}
where the dot denotes inner product;
then the random-walk Metropolis algorithm with $q$
bounded away from zero on a neighborhood of zero is geometrically ergodic.
\end{thmm}
We say \emph{$\pi$ satisfies the curvature condition}
to mean \eqref{condjh-curvature} holds.
This means the contours of $\pi$ are approximately locally linear near
infinity.

Theorem~\ref{thmjh-main}, although useful, covers neither exponentially
light densities, which arise in Bayesian categorical data analysis with
canonical parameters and conjugate priors
(Section~\ref{secbounded-exponential}),
nor sub-exponentially light densities, which
arise in Bayesian analysis of Cauchy location models using flat improper
priors on the location parameters (Section~\ref{seccauchy-location}).
\citet{roberts-tweedie} do cover exponentially light densities, but their
theorems are very difficult to apply [\citet{jarner-hansen} show that
\citet{roberts-tweedie} incorrectly applied their own theorem in one
case].

The key idea of this paper is to use the change-of-variable theorem in
conjunction with Theorem~\ref{thmjh-main} to get results that
Theorem~\ref{thmjh-main} does not give directly.
Suppose~$\pi_\beta$ is the (possibly multivariate) target density of interest.
We instead simulate a Markov chain having invariant density
%
%
%e4 #&#
\begin{equation}
\label{eqchange} \pi_\gamma(\gamma) = \pi_\beta \bigl(h(\gamma)
\bigr) \bigl\vert\det\nabla h(\gamma) \bigr\vert,\vadjust{\goodbreak}
\end{equation}
where $h$ is a diffeomorphism. If $\pi_\beta$ is the density of
the random vector $\beta$, then~$\pi_\gamma$ is the density of
the random vector $\gamma= h^{- 1}(\beta)$.
We find conditions
on the transformation $h$ that make $\pi_\gamma$ super-exponentially
light and satisfy the curvature condition.
Then by Theorem~\ref{thmjh-main}, the simulated Markov chain
$\gamma_1$, $\gamma_2$, $\ldots$ is geometrically ergodic.
It is easy to see (Appendix~\ref{secisomorphic}) that
the Markov chain $\beta_i = h(\gamma_i)$, $i = 1, 2, \ldots,$
is also geometrically ergodic. Thus we achieve geometric ergodicity
indirectly, doing a change-of-variable yielding a density that by
Theorem~\ref{thmjh-main} has a geometrically ergodic random-walk
Metropolis sampler, sampling that distribution, and then using the
inverse change-of-variable to get back to the variable of interest.

This indirect procedure has no virtues other than that Metropolis
random-walk samplers are well-understood and user-friendly
and that we have Theorem~\ref{thmjh-main} to build on.
There is other literature using drift conditions to prove geometric ergodicity
of Markov chain samplers
[\citet{geyer-moller},
\citet{rosenthal-james-stein},
\citet{hobert-geyer},
\citet{jones-hobert},
\citet{roy-hobert},
\citet{tan-hobert},
\citet{johnson-jones}]
but for Gibbs samplers or
other samplers for specific statistical models, hence not having the wide
applicability of random-walk Metropolis samplers.
There is also other literature about using variable transformation to improve
the convergence properties of Markov chain samplers
[\citet{roberts-sahu},
\citet{papaspiliopoulos-roberts-skold},
\citet{papaspiliopoulos-roberts}]
but for Gibbs samplers not having the wide applicability of random-walk
Metropolis samplers.

It is important to understand that the necessary condition mentioned above
[\citet{jarner-tweedie}] places a limit on what can be done without variable
transformation. If $\pi_\beta$ does not have a moment generating function
(any Student $t$ distribution, e.g.), then no random-walk Metropolis
sampler for it can be geometrically ergodic (no matter what proposal
distribution is used). Thus if we use a random-walk Metropolis sampler,
then we must also use variable transformation to obtain geometric ergodicity.

We call a function $h \dvtx\real^k \to\real^k$ \emph{isotropic} if
it has
the form
%
%
%e5 #&#
\begin{equation}
\label{eqh} h(\gamma) = \cases{\displaystyle f\bigl(\vert\gamma \vert\bigr) \frac
{\gamma
}{\vert\gamma \vert}, &\quad $
\gamma \neq0$, \vspace*{2pt}
\cr
0, &\quad $\gamma= 0$}
\end{equation}
for some function $f \dvtx(0, \infty) \to(0, \infty)$.
To simplify the theory, we restrict attention to
$h$ that are isotropic diffeomorphisms, meaning $h$ and $h^{- 1}$
are both continuously differentiable, having the further property that
$\det(\nabla h)$ and $\det(\nabla h^{- 1})$ are also continuously
differentiable.

As with the restriction to
$\pi$ that are strictly positive and continuously differentiable
used by \citet{mengersen-tweedie}, \citet{roberts-tweedie}
and \citet{jarner-hansen}, this restriction is arbitrary.
It is not necessary to achieve geometric ergodicity; it merely simplifies
proofs.
However, the proofs are already\vadjust{\goodbreak} very complicated even with
these two restrictions.
Although both these restrictions could be relaxed, that would make
the proofs even more complicated.
Since many applications can be fit into our framework, perhaps after
a change-of-variable to yield $\pi_\beta$ that is
strictly positive and continuously differentiable,
we choose to not complicate our proofs further.

Isotropic transformations \eqref{eqh} shrink toward or expand away from
the origin of the state space. In practice, they should be combined with
translations so they can shrink toward or expand away from arbitrary points.
Since translations induce isomorphic Markov chains
(Appendix~\ref{secisomorphic}), they do not affect
the geometric ergodicity properties of random-walk Metropolis samplers.
Hence we ignore them until Section~\ref{secdiscuss}.

Our variable-transformation method is easily implemented using
the R package \texttt{mcmc} [\citet{mcmc-R-package}] because that
package simulates Markov chains having equilibrium
density $\pi$ specified by a user-written function, which can incorporate
a variable transformation, and outputs an arbitrary functional of the
Markov chain specified by another user-written function, which can incorporate
the inverse transformation.

A referee pointed out that one can think of our transformation method
differently: as describing a Metropolis--Hastings
algorithm in the original parameterization.
This seems to avoid variable transformation but does not,
because its proposals have the form
$h (h^{- 1}(\beta) + z )$, where $\beta$ is the current state,
and $z$ is a simulation from the Metropolis $q$.
This uses $h$ and $h^{- 1}$ in every iteration,
whereas the scheme we describe uses only $h$ to run the Markov chain
for $\gamma$ and to map it back to $\beta$, needing $h^{- 1}$ only once
to determine the inital state $\gamma_1 = h^{- 1}(\beta_1)$ of the Markov
chain. Nevertheless, it is of some
theoretical interest that this provides hitherto unnoticed examples
of geometrically ergodic Metropolis--Hastings algorithms.

%s2 #&#
\section{Variable transformation}

%s2.1 #&#
\subsection{Positivity and continuous differentiability} \label{secpos-c1}

For the change-of-variable \eqref{eqchange} we need to know when
the transformed density $\pi_\gamma$ is positive and continuously
differentiable assuming the original density $\pi_\beta$ has these
properties. If $h$ is a diffeomorphism, then the first term
on the right-hand side will be continuously differentiable by the chain rule.
Since $\nabla h^{- 1}$ is the matrix inverse of $\nabla h$ by the inverse
function theorem, $\det(\nabla h)$ can never be zero. Hence $h$ being
a diffeomorphism is enough to imply positivity of $\pi_\gamma$.

Since $\det(A)$ is
continuous in $A$, being a polynomial function of the components of $A$,
$\det(\nabla h)$ can never change sign. We restrict attention to $h$
such that $\det(\nabla h)$ is always positive, so the absolute value in
\eqref{eqchange} is unnecessary. Then we have
%
%
%e6 #&#
%e7 #&#
\begin{eqnarray}
\label{eqdensity-trans-log} \log\pi_\gamma(\gamma) & =& \log
\pi_\beta \bigl( h (\gamma) \bigr) + \log\det \bigl( \nabla h ( \gamma)
\bigr),
\\
\label{eqdensity-trans-dee-log} \nabla\log\pi_\gamma(\gamma) & = &\nabla
( \log\pi_\beta) \bigl( h(\gamma) \bigr) \nabla h(\gamma) + \nabla\log
\det \bigl( \nabla h ( \gamma) \bigr).\vadjust{\goodbreak}
\end{eqnarray}
It is clear from \eqref{eqdensity-trans-dee-log} that $\log\pi_\gamma$,
and hence $\pi_\gamma$ is continuously differentiable if $h$ is
a diffeomorphism, and $\det(\nabla h)$ is continuously differentiable.

%s2.2 #&#
\subsection{Isotropic functions}

In the transformation method, the induced density, $\pi_\gamma$ will need
to satisfy the smoothness conditions of Theorem~\ref{thmjh-main}. We
require the original density, $\pi_\beta$ to satisfy the smoothness
conditions of Theorem~\ref{thmjh-main}. The smoothness conditions will be
satisfied for $\pi_\gamma$ if the isotropic transformations are
diffeomorphisms with continuously differentiable Jacobians. The
assumptions of the following lemma provide conditions on isotropic
functions to guarantee that $\pi_\gamma$ is positive and continuously
differentiable whenever $\pi_\beta$ is.

%
%le1 #&#
\begin{lem}\label{lemf-conds}
Let $h\dvtx\real^k\rightarrow\real^k$ be an isotropic function
given by \eqref{eqh} with $f \dvtx[0, \infty) \rightarrow[0, \infty)$
invertible and continuously differentiable with one-sided derivative
at zero such that
%
%
%e8 #&#
\begin{equation}
\label{condf-first} f'(s) > 0, \qquad s \ge0.
\end{equation}
Then
%
%
%e9 #&#
\begin{equation}
\label{eqgamma-ratio} \frac{\gamma}{\vert\gamma \vert} = \frac
{h(\gamma)}{\vert h(\gamma) \vert},\qquad  \gamma\neq0,
\end{equation}
$f$ is a diffeomorphism, $h$ is a diffeomorphism and
%
%
%e10 #&#
\begin{equation}
\label{eqh-inverse} h^{- 1}(\beta) = \cases{\displaystyle f^{- 1}\bigl(\vert
\beta \vert\bigr) \frac{\beta}{\vert\beta \vert}, & \quad$\beta\neq0$, \vspace*{2pt}
\cr
0, &\quad $ \beta= 0$}
\end{equation}
and
%
%
%e11 #&#
\begin{equation}
\label{eqdee-h-nonzero} \nabla h(\gamma) = \frac{f(\vert\gamma
\vert)
\I_k}{\vert\gamma \vert} + \biggl[
f^\prime\bigl(\vert\gamma \vert\bigr) - \frac{f(\vert\gamma \vert)}{\vert
\gamma \vert} \biggr]
\frac{\gamma\gamma^T}{\vert\gamma \vert^2},\qquad  \gamma\neq0,
\end{equation}
where $\I_k$ is the $k\times k$ identity matrix, and
%
%
%e12 #&#
\begin{equation}
\label{eqdee-h-zero} \nabla h(0) = f'(0) \I_k.
\end{equation}
Moreover
%
%
%e13 #&#
\begin{equation}
\label{eqdet-dee-h} \det \bigl(\nabla h(\gamma) \bigr) = \cases{
\displaystyle f'\bigl(\vert\gamma \vert\bigr) \biggl( \frac{f(\vert\gamma \vert)}{\vert
\gamma \vert}
\biggr)^{k-1}, &\quad $\gamma\neq0,$ \vspace*{2pt}
\cr
f'(0)^{k},
& \quad $\gamma= 0$}
\end{equation}
and, under the additional assumption that $f$ is twice continuously
differentiable with one-sided derivatives at zero and
%
%
%e14 #&#
\begin{equation}
\label{condf-second} f''(0) = 0,
\end{equation}
\eqref{eqdet-dee-h} is continuously differentiable.
\end{lem}
The proof of this lemma is in Appendix~\ref{secproof-lemma-isotropic}.

%s2.3 #&#
\subsection{Inducing lighter tails} \label{seclighter-tails}

Define $f \dvtx[0, \infty) \to[0, \infty)$ by
%
%
%e15 #&#
\begin{equation}
\label{eqf1} f(x) = \cases{ x, &\quad  $x < R$, \vspace*{2pt}
\cr
x + (x -
R)^p, & \quad $x \geq R$,}
\end{equation}
where $R \geq0$ and $p > 2$.
It is clear that \eqref{eqf1} satisfies the assumptions
of Lemma~\ref{lemf-conds}.
% changed to make parallel construction with comment after Theorem 3.

%
%th2 #&#
\begin{thmm}\label{thmexponential-to-super}
Let $\pi_\beta$ be an exponentially light density on $\real^k$, and
let $h$ be defined by \eqref{eqh} and \eqref{eqf1}. Then
$\pi_\gamma$ defined by \eqref{eqchange}
is super-exponentially light.
\end{thmm}

Proof of Theorem~\ref{thmexponential-to-super} is in Appendix
\ref{secproof-induced}.

Now define $f \dvtx[0, \infty) \to[0, \infty)$ by
%
%
%e16 #&#
\begin{equation}
\label{eqf2} f(x) = \cases{\displaystyle e^{bx} - \frac{e}{3}, &\quad  $\displaystyle x >
\frac{1}{b}$, \vspace*{2pt}
\cr
\displaystyle x^3 \frac{b^3e}{6} + x
\frac{be}{2}, &\quad  $\displaystyle x \leq\frac{1}{b},$}
\end{equation}
where $b > 0$.
It is clear that \eqref{eqf2} satisfies the assumptions
of Lemma~\ref{lemf-conds}.

%
%th3 #&#
\begin{thmm}\label{thmsub-exponential-to-exponential}
Let $\pi_\beta$ be a sub-exponentially light density on $\real^k$, and
suppose there exist $\alpha> k$ and $R < \infty$ such that
%
%
%e17 #&#
\begin{equation}
\label{condsub-bounded} \frac{\beta}{\vert\beta \vert} \cdot \nabla \log\pi_\beta(
\beta) \leq- \frac{\alpha}{\vert\beta \vert}, \qquad \vert\beta \vert > R.
\end{equation}
Let $h$ be defined by \eqref{eqh} and \eqref{eqf2}. Then
$\pi_\gamma$ defined by \eqref{eqchange}
is exponentially light.
\end{thmm}

Proof of Theorem~\ref{thmsub-exponential-to-exponential} is in Appendix
\ref{secproof-induced}.

Condition \eqref{condsub-bounded} is close to sharp.
For example, if $\pi_\beta$ looks like a multivariate $t$ distribution
%
%
%e18 #&#
\begin{equation}
\label{eqfubar} \pi_\beta(t) = \bigl[ 1 + (t - \mu)^T
\Sigma^{- 1} (t - \mu) \bigr]^{- (v + k) / 2}
\end{equation}
[compare with \eqref{eqt-multivariate} in Section~\ref{sect-multivariate}],
then \eqref{condsub-bounded} holds with $\alpha= k + v$, and
\eqref{eqfubar} is integrable if and only if $v > 0$.

Moreover, an exponential-type isotropic transformation like \eqref{eqf2}
is necessary to obtain a super-exponentially light $\pi_\gamma$ when
$\pi_\beta$ is a multivariate $t$ distribution.
Direct calculation shows that no polynomial-type isotropic transformation
like \eqref{eqf1} does the job.

%
%co1 #&#
\begin{cor}\label{corcomposition}
Let $\pi_\beta$ satisfy the conditions of
Theorem~\ref{thmsub-exponential-to-exponential}, and let $h$ be defined
as the composition of those used in
Theorems~\ref{thmexponential-to-super}
and~\ref{thmsub-exponential-to-exponential}; that is, if we denote the~$h$
used in Theorem~\ref{thmexponential-to-super} by $h_1$ and denote
the~$h$ used in Theorem~\ref{thmsub-exponential-to-exponential} by~$h_2$,
then in this corollary we are using $h = h_2 \circ h_1$ and the
change of variable is $\gamma= h_1^{-1}(h_2^{-1}(\beta))$. Then
$\pi_\gamma$ defined by \eqref{eqchange} is super-exponentially light.
\end{cor}
\begin{pf}
The proof follows directly from Theorems~\ref{thmexponential-to-super}
and~\ref{thmsub-exponential-to-exponential}.
\end{pf}

%s2.4 #&#
\subsection{Curvature conditions}\label{seccurvature}

As seen in \citet{jarner-hansen}, Example~5.4, being super-exponentially
light is \emph{not} a sufficient condition for the geometric
ergodicity of
a random-walk Metropolis algorithm. \citet{jarner-hansen} provide
sufficient conditions for super-exponentially light densities. In this
section, we provide sufficient conditions for sub-exponentially light and
exponentially light densities, such that, using the transformations from
Section~\ref{seclighter-tails} the induced super-exponential densities
will satisfy the \citet{jarner-hansen} sufficient conditions.

%
%th4 #&#
\begin{thmm}\label{thmcurvature-exponential}
Let $\pi_\beta$ be an exponentially light density on $\real^k$, and
suppose that $\pi_\beta$ satisfies either of the following conditions:
\begin{longlist}[(ii)]%[label={\normalfont(\roman{*})}]
\item[(i)] $\pi_\beta$ satisfies the
curvature condition
\eqref{condjh-curvature}, or
\item[(ii)] $\vert\nabla\log \pi_\beta
(\beta) \vert$ is bounded as $\vert\beta \vert$ goes to infinity.
\end{longlist}
Let $h$ be defined by \eqref{eqh} and \eqref{eqf1}. Then
$\pi_\gamma$ defined by \eqref{eqchange}
satisfies the curvature condition
\eqref{condjh-curvature}.
\end{thmm}

Proof of Theorem~\ref{thmcurvature-exponential} is in Appendix
\ref{secproof-curvature}.

For exponentially light $\pi_\beta$,
condition~(ii)
implies condition~(i). In practice,
condition~(ii) may be easier to check than
condition~(i)
(as in Section~\ref{secbounded-exponential}).

%
%th5 #&#
\begin{thmm}\label{thmcurvature-sub-exponential}
Let $\pi_\beta$ be a sub-exponentially light density on $\real^k$, and
suppose there exist $\alpha> k$ and $R < \infty$ such that
%
%
%e19 #&#
\begin{equation}
\label{condsub-magnitude-bdd} \bigl\vert\nabla\log\pi_\beta(\beta) \bigr\vert
\leq \frac{\alpha}{\vert\beta \vert},\qquad  \vert\beta \vert > R.
\end{equation}
Let $h$ be defined by \eqref{eqh} and \eqref{eqf2}. Then
$\pi_\gamma$ defined by \eqref{eqchange} satisfies
condition \textup{(ii)} of Theorem
\ref{thmcurvature-exponential} with $\beta$ replaced by $\gamma$.
\end{thmm}

Proof of Theorem~\ref{thmcurvature-sub-exponential} is in Appendix
\ref{secproof-curvature}.

Condition \eqref{condsub-magnitude-bdd}, like \eqref{condsub-bounded},
is close to sharp. If $\pi_\beta$ has the form \eqref{eqfubar},
then \eqref{condsub-magnitude-bdd} holds with $\alpha= k + v$, and
\eqref{eqfubar} is integrable if and only if $v > 0$.

%
%co2 #&#
\begin{cor}\label{corcurvature-sub-to-super}
Let $\pi_\beta$ satisfy the conditions of
Theorems~\ref{thmsub-exponential-to-exponential} and~\ref
{thmcurvature-sub-exponential},
and let $h$ be defined as the composition of those used in
Theorems~\ref{thmcurvature-exponential}
and~\ref{thmcurvature-sub-exponential}, that is, if we denote the $h$
used in Theorem~\ref{thmcurvature-exponential} by $h_1$ and denote the
$h$ used in Theorem~\ref{thmcurvature-sub-exponential} by $h_2$, then in
this corollary we are using $h = h_2 \circ h_1$ and the change of
variable is $\gamma= h_1^{-1}(h_2^{-1}(\beta))$. Then $\pi_\gamma$
defined by \eqref{eqchange} satisfies the curvature condition
\eqref{condjh-curvature}.
\end{cor}
\begin{pf}
This follows directly from Theorems~\ref{thmcurvature-sub-exponential}
and~\ref{thmcurvature-exponential}.
\end{pf}

To verify that a variable transformation \eqref{eqh} produces geometric
ergodicity, one uses Theorems~\ref{thmexponential-to-super}
and~\ref{thmcurvature-exponential} when the given target density $\pi_\beta$
is exponentially light.
To verify that a variable transformation \eqref{eqh} produces geometric
ergodicity, one uses Corollaries~\ref{corcomposition}
and~\ref{corcurvature-sub-to-super} when the given target density $\pi_\beta$
is sub-exponentially light.
(When the given target density $\pi_\beta$ is super-exponentially light
one does not need variable transformation to obtain geometric ergodicity
if $\pi_\beta$ also satisfies the curvature condition.)\

% Rule 1. Only one or two readers, those who are going to extend the
%theory
% theory in this paper in some way, are ever going to read anything but
%the
% introduction and the discussion. Therefore put every point you want
%readers
% to understand in the introduction or discussion, including stuff like
%below
% (why we chose the strategy we chose). You should consider yourself
%lucky
% if a reader reads the introduction and discussion. Most readers will
%read
% only the title and abstract (or even just the title) and skip the
%rest.
% Think about someone seeing our title in Google Scholar and skipping
%the
% link. That's the vast majority of "readers".
%
% Thinking about that, I changed the title.
%
% \textbf{VERY ROUGH, REWORD}
%
% Again, as in Theorem~\ref{thmcurvature-sub-exponential} we are not
% directly showing that the transformation induces a density that
%satisfies
% the conditions of Theorem~\ref{thmjh-main}. Instead, because we are
% composing multiple transformations, we are showing that the
%transformation
% will induce a density that satisfies the conditions of Theorem
%~\ref{thmcurvature-exponential}, which shows that the density induced
%by
% the composition of transformations will satisfy the conditions of
%Theorem
%~\ref{thmjh-main}.

%s3 #&#
\section{Examples}

% I think we want subsections

%s3.1 #&#
\subsection{Exponential families and conjugate priors}
\label{secbounded-exponential}

In this section we study Bayesian inference for exponential families using
conjugate priors, in particular, the case where
the natural statistic is bounded in some direction, and the natural
parameter space is all of $\real^k$.
Examples include logistic
regression, Poisson regression with log link function
and log-linear models in categorical data analysis.
In this case, we find that the posterior density, when it exists,
is exponentially light
and satisfies the curvature condition. Hence variable transformation using
\eqref{eqh} and~\eqref{eqf1} makes the random-walk Metropolis sampler
geometrically ergodic.

An exponential family is a statistical model having log likelihood of the
form
\[
y \cdot\beta- c(\beta),
\]
where the dot denotes inner product, $y$ is a vector statistic, $\beta
$ is
a vector parameter and the function $c$ is called the cumulant function
of the family. A statistic $y$ and parameter $\beta$ that give a log likelihood
of this form are called \emph{natural} or \emph{canonical}.
If $y_1, \ldots, y_n$ are independent and identically
distributed observations from the family and $\bar{y}_n$ their average,
then the log likelihood for the sample of size $n$ is
\[
n \bar{y}_n \cdot\beta- n c(\beta).
\]
The log unnormalized posterior when using conjugate priors is
%
%
%e20 #&#
\begin{equation}
\label{eqlup} w(\beta) = (n \bar{y}_n + \nu\eta) \cdot\beta- (n +
\nu) c(\beta),
\end{equation}
where $\nu$ is a scalar hyperparameter, and $\eta$ is a vector hyperparameter
[\citet{diaconis-ylvisaker}, Section~2]. When simulating the posterior
using MCMC, the unnormalized density of the target distribution is
$\pi(\beta) = e^{w(\beta)}$.

The \emph{convex support} of an exponential family is
the smallest closed convex set containing the natural statistic
with probability one. (This does not depend on which distribution in the
exponential family we use because they are all mutually absolutely continuous.)
Theorem~1 in \citet{diaconis-ylvisaker} says that the posterior exists;
that is, $e^{w(\beta)}$ is integrable,
where $w(\beta)$ is given by \eqref{eqlup},
if and only if $n + \nu> 0$ and $(n \bar{y}_n + \nu\eta) / (n + \nu
)$ is
an interior point of the convex support.
(Of course, this always happens when using a proper prior, i.e.,
when $\nu> 0$ and $\eta/ \nu$ is an interior point of the convex support.)

Theorem~9.13 in \citet{barndorff} says that this same condition holds
if and only if the log unnormalized posterior \eqref{eqlup}
achieves its maximum at a unique point, the posterior mode,
call it $\tilde{\beta}_n$.
(Ostensibly, this theorem applies only to log likelihoods
of exponential families not to log unnormalized posteriors with conjugate
priors, but since the latter have the same algebraic form as the former,
it actually does apply to the latter.)

From the properties of exponential families [\citet{barndorff}, Theorem~8.1],
%
%
%e21 #&#
\begin{equation}
\nabla c(\beta)  = E_\beta(Y). \label{eqbartlett-1} % \\
% \nabla^2 c(\beta) & = \var_\beta(Y)
% \label{eqbartlett-2}
\end{equation}
It follows that
%
%
%e22 #&#
\begin{equation}
\label{eqlup-deriv} \nabla\log\pi(\beta) = \nabla w(\beta) = n
\bar{y}_n + \nu\eta- (n + \nu) E_\beta(Y).
\end{equation}
Suppose that the natural statistic is bounded in some direction, that is,
there exists a nonzero vector $\delta$ and real number $b$ such that
$y \cdot\delta\le b$ for all $y$ in the convex support.
It follows that $E_\beta(Y) \cdot\delta\le b$.
Then
\begin{eqnarray*}
\limsup_{\vert\beta \vert \rightarrow\infty} \frac{\beta}{\vert
\beta \vert} \cdot\nabla\log\pi(\beta) & \ge&
\limsup_{s \rightarrow\infty} \frac{s \delta}{\vert s \delta
\vert} \cdot \bigl[ n \bar{y}_n +
\nu\eta- (n + \nu) E_{s \delta}(Y) \bigr]
\\
& \ge&\frac{(n \bar{y}_n + \nu\eta) \cdot\delta- (n + \nu)
b}{\vert\delta \vert}.
\end{eqnarray*}
Hence \eqref{eqtails} is not $- \infty$ and the target distribution
is not super-exponentially light.

When the convex support has nonempty interior,
the cumulant function $c$ is strictly convex [\citet{barndorff}, Theorem~7.1].
Hence \eqref{eqlup} is a strictly concave function.
It follows from this that $\nabla c$ is a strictly multivariate monotone
function, that is,
%
%
%e23 #&#
\begin{equation}
\label{eqmono-general} \bigl[ \nabla c(\beta_1) - \nabla c(
\beta_2) \bigr] \cdot(\beta_1 - \beta_2) > 0,\qquad
\beta_1 \neq\beta_2
\end{equation}
[\citet{rockafellar-wets}, Theorem~2.14 and Chapter~12].
It follows that
%
%
%e24 #&#
\begin{equation}
\label{eqmono} \nabla w(\beta) \cdot\frac{\beta- \tilde{\beta}_n}{
\vert\beta- \tilde{\beta}_n \vert} < 0,\qquad  \beta\neq\tilde{
\beta}_n,
\end{equation}
where $w$ is given by \eqref{eqlup}, because $\nabla w(\tilde{\beta}_n)
= 0$.
Let $B$ denote the boundary and~$E$ denote the exterior of the ball of unit
radius centered at $\tilde{\beta}_n$.
Since $c$ is infinitely differentiable
[\citet{barndorff}, Theorem~7.2],\vadjust{\goodbreak} so is $w$, and the left-hand side
of \eqref{eqmono} is a continuous function of $\beta$. Since $B$ is compact,
the left-hand side of \eqref{eqmono} achieves its maximum over $B$, which
must be negative, say $- \varepsilon$.
For any $\beta\in E$ we have $t \beta+ (1 - t) \tilde{\beta}_n \in B$
when $t = 1 / \vert\beta- \tilde{\beta}_n \vert$.
By \eqref{eqmono-general} we have
\[
\bigl[ \nabla w(\beta) - \nabla w \bigl( t \beta+ (1 - t) \tilde{
\beta}_n \bigr) \bigr] \cdot\frac{\beta- \tilde{\beta}_n}{
\vert\beta- \tilde{\beta}_n \vert} < 0
\]
%
% deleted \beta\in E because t depends on \beta
because
\[
\beta- \bigl[ t \beta+ (1 - t) \tilde{\beta}_n \bigr] = (1 - t) (
\beta- \tilde{\beta}_n )
\]
is parallel to $\beta- \tilde{\beta}_n$. Thus
\[
\nabla w(\beta) \cdot\frac{\beta- \tilde{\beta}_n}{ \vert\beta
- \tilde{\beta}_n \vert} < - \varepsilon, \qquad \beta\in E
\]
and
\[
\limsup_{\beta\to\infty} \nabla w(\beta) \cdot\frac{\beta-
\tilde{\beta}_n}{ \vert\beta- \tilde{\beta}_n \vert} \le- \varepsilon,
\]
and this is easily seen to be equivalent to the unnormalized density
\eqref{eqlup} being exponentially light.

Now we check the curvature condition \eqref{condjh-curvature} for exponential
families. In case the natural statistic is bounded in all directions, as
in logistic regression and log-linear models, the curvature condition follows
directly because the family satisfies condition~(ii)
of Theorem~\ref{thmcurvature-exponential} because $\nabla\log\pi
(\beta)$
is \eqref{eqlup-deriv}, and this is bounded.
In case the natural statistic is bounded in some directions but not all
directions, as in Poisson regression, we have to work harder and use
condition~(i)
of Theorem~\ref{thmcurvature-exponential}. Because
\[
\nabla\log\pi(\beta) = \frac{\nabla\pi(\beta)}{\pi(\beta)},
\]
we have
\[
\frac{\nabla\pi(\beta)}{\vert\nabla\pi(\beta) \vert} = \frac
{\nabla w(\beta)}{\vert\nabla w(\beta) \vert},
\]
where $\nabla w(\beta)$ is given by \eqref{eqlup-deriv}.
And from \eqref{eqmono} and $\nabla w(\beta) \neq0$
for $\beta\neq\tilde{\beta}_n$, we obtain
%
%
%e25 #&#
\begin{equation}
\label{eqmono-curvature} \frac{\nabla w(\beta)}{\vert\nabla
w(\beta) \vert} \cdot\frac{\beta- \tilde{\beta}_n}{
\vert\beta- \tilde{\beta}_n \vert} < 0,\qquad  \beta
\neq\tilde{\beta}_n,
\end{equation}
and the rest of the proof that $\pi$ satisfies the curvature condition
is just like the proof that it is exponentially light given above except
that \eqref{eqmono-curvature} replaces \eqref{eqmono}.

%s3.2 #&#
\subsection{Multinomial logit regresion with a conjugate prior}
\label{secmultinomial-logit}
This example is a special case of the example in
Section~\ref{secbounded-exponential}.\vadjust{\goodbreak}

In multinomial logit regression, using a conjugate prior is equivalent
to adding prior counts to the data cells.
For observations $1,\dots,L$, represent these
prior counts as $\xi_l \nu_l$ where $\xi_l$ is a vector giving the prior
probability for each response for the $l$th observation, and $\nu_l$
is the prior sample size. For the $l$th observation,
let the vector $Y^l$ represent the counts in each response category,
$N^l = \sum_i Y^l_i$ be the sample size and $M^l$ be the model matrix.
The log unnormalized posterior density for the regression parameter
$\beta$
is given by
%
%
%e26 #&#
\begin{equation}
\label{eqmult-post-conj-complicated} \qquad \pi(\beta|y,n,\xi,\nu) \propto\exp \Biggl\{
\sum_{l=1}^L \bigl(y^l +
\xi^l \nu^l \bigr) \cdot M^l \beta-
\bigl(n^l + \nu^l \bigr) \log \biggl( \sum
_j e^{M_{j\cdot}^l \beta} \biggr) \Biggr\},
\end{equation}
where $M_{j\cdot}^l$ is the $j$th row of the matrix $M^l$.
So long as $y^l_i + \xi^l_i \nu^l$ is positive for all $i$ and
$l$---there is data (actual plus prior) in all cells---$\pi$
will be
exponentially light, and satisfy condition~\eqref{condjh-curvature}.
Hence a random-walk Metropolis algorithm for the density induced by the
approach in Theorems~\ref{thmexponential-to-super}
and~\ref{thmcurvature-exponential} will be geometrically ergodic.

%s3.3 #&#
\subsection{Multivariate T distributions}
\label{sect-multivariate}

The density of a multivariate $t$ distribution on $\real^k$ with $v$
degrees of freedom, location parameter vector $\mu$ and scale parameter
matrix $\Sigma$ is given by
%
%
%e27 #&#
\begin{equation}
\label{eqt-multivariate}\qquad \pi_\beta(t) = \frac{\Gamma[(v+k)/2
]}{\Gamma[v/2 ]
(v\pi)^{k/2} \det(\Sigma)^{1/2}} \biggl[ 1 +
\dfrac{1}{v} (t - \mu)^T \Sigma^{-1}(t - \mu)
\biggr]^{-(v+k) / 2}
\end{equation}
so
%
%
%e28 #&#
\begin{equation}
\label{eqgosset} \nabla\log\pi_\beta(t) = \frac{-(v+k) \Sigma^{-1}(t - \mu)} {
v + (t - \mu)^T \Sigma^{-1}(t - \mu)},
\end{equation}
which implies
%
%
%e29 #&#
\begin{equation}
\label{eqgosset-too} t \cdot\nabla\log\pi_\beta(t) \to- (v + k),\qquad
\mbox{as $t \to\infty$},
\end{equation}
so \eqref{eqt-multivariate} is sub-exponentially light.
% sentence moved below

The condition of Theorem~\ref{thmsub-exponential-to-exponential}
is also implied by \eqref{eqgosset-too}. To check the condition of
Theorem~\ref{thmcurvature-sub-exponential} we calculate
\[
\bigl\vert\nabla\log\pi_\beta(t) \bigr\vert^2 \le\frac{(v + k)^2 \lambda_{\mathrm{max}}^2 \vert t - \mu \vert^2} {
(\lambda_{\mathrm{min}} \vert t - \mu \vert^2)^2},
\]
where $\lambda_{\mathrm{max}}$ and $\lambda_{\mathrm{min}}$ are
the largest and smallest eigenvalues of $\Sigma^{- 1}$. Hence
\[
\bigl\vert\nabla\log\pi_\beta(t) \bigr\vert \le\frac{(v + k) \lambda_{\mathrm
{max}}}{\lambda_{\mathrm{min}} \vert t - \mu \vert},
\]
and the condition of Theorem~\ref{thmcurvature-sub-exponential} also holds.
So a random-walk Metropolis algorithm for
the induced density $\pi_\gamma$ that uses the transformation
described in
Corollaries~\ref{corcomposition} and~\ref{corcurvature-sub-to-super}
will be geometrically ergodic, and the inverse transformed Markov chain
will be geometrically ergodic for $\pi_\beta$.
Since the multivariate $t$ distribution does not have a moment generating
function, no random-walk Metropolis algorithm for $\pi_\beta$ is
geometrically ergodic [\citet{jarner-tweedie}].
Variable transformation is essential.

The case $k = 1$ gives the univariate $t$ distribution, which has been
widely used as an example of a Harris ergodic random-walk Metropolis algorithm
that is not geometrically ergodic
[\citet{mengersen-tweedie},
\citet{jarner-hansen},
\citet{jarner-tweedie},
\citet{jarner-roberts}].

%s3.4 #&#
\subsection{Cauchy location models and flat priors}
\label{seccauchy-location}

The $t$ distribution with one degree of freedom is the Cauchy distribution.
Consider a Cauchy location family with flat prior, so the posterior density
for sample size one is again a Cauchy distribution
\[
\pi_\beta(\mu) = \frac{1}{\pi} \cdot\frac{1}{1 + (x - \mu)^2},
\]
and, this being a special case of the preceding section, this density
is sub-exponentially light.

For a sample of size $n$ the unnormalized posterior density is
\[
\pi_\beta(\mu) = \prod_{i = 1}^n
\frac{1}{1 + (x_i - \mu)^2}
\]
and the posterior distribution is no longer a brand name distribution.
It is still easily shown to be sub-exponentially light and to satisfy
the conditions of Theorems~\ref{thmsub-exponential-to-exponential}
and~\ref{thmcurvature-sub-exponential}.

%s4 #&#
\section{Discussion} \label{secdiscuss}

The transformations in
Theorems~\ref{thmexponential-to-super} and~\ref
{thmsub-exponential-to-exponential}
will always induce a density with tails at least as light as the original
density. If the original density satisfies the curvature condition, then
the transformation using the transformation from
Theorem~\ref{thmexponential-to-super} will induce a density that satisfies
the curvature condition. Thus applying the transformation from
Theorem~\ref{thmexponential-to-super} to a super-exponentially light
density that satisfies the curvature condition will induce another
super-exponentially light density that satisfies the curvature condition.
We do not recommend transformation when the original density already
satisfies the conditions of Theorem~\ref{thmjh-main},
but it seems this will do no harm.

The transformation method introduced here can be mixed blessing.
It can produce geometric ergodicity, but may cause other problems. For
example, $\pi_\gamma$ given by \eqref{eqchange} can be multimodal when
$\pi_\beta$ is unimodal. Thus we want a less extreme member of the family
of transformations that does the job. The idea is to pull in the tails
enough to get geometric ergodicity without much affecting the main part of
the distribution. Although very extreme transformations work in theory,
they are problematic in practice due to inexactness of computer arithmetic.

As mentioned in the \hyperref[sec1]{Introduction}, in practice one
combines the
transformations introduced in Section~\ref{seclighter-tails}
with translations. Let $t_\lambda$ denote the translation
$x \mapsto x + \lambda$.
Then in the exponentially light $\pi_\beta$ case, we use the transformation
$h = t_\lambda\circ h_{R, p}$, where $h_{R, p}$ is the $h$ defined by
\eqref{eqh} and \eqref{eqf1}, so the change-of-variable is
$\gamma= h_{R, p}^{- 1}(\beta- \lambda)$.
This gives users three adjustable constants,
$\lambda$, $R$ and $p$, to experiment with to improve the mixing
of the sampler. If $\pi_\beta$ satisfies the assumptions of
Theorems~\ref{thmexponential-to-super} and~\ref{thmcurvature-exponential},
then any valid values of $\lambda$, $R$ and $p$ result in a geometrically
ergodic sampler. Observe that the restriction of this $h$ to the ball of
radius $R$ centered at $\lambda$ is a translation, which does not
affect the
shape of the distribution. Thus one wants to choose $\lambda$ near the
center of the distribution (perhaps the mode of $\pi_\beta$, if it
has one)
and $R$ large enough so that a large part of the probability is in this ball
where the shape is unchanged. The parameter $p$ should always be chosen to
be small, say 3 or 2.5 (recall $p > 2$ is required), 3 is a good
choice as then $f$ has a closed-form expression for its inverse.

In the sub-exponentially
light $\pi_\beta$ case, we use the transformation
$h = t_\lambda\circ h_b \circ h_{R, p}$, where $h_b$ is the $h$
defined by
\eqref{eqh} and \eqref{eqf2}, and the other two transformations are
as above,
so the change-of-variable is
$\gamma= h_{R, p}^{- 1} (h_b^{- 1}(\beta- \lambda) )$.
This gives users four adjustable constants,
$\lambda$, $R$, $p$ and $b$ to experiment with to improve the mixing
of the sampler. If $\pi_\beta$ satisfies the assumptions of
Corollaries~\ref{corcomposition} and~\ref{corcurvature-sub-to-super},
then any valid values of $\lambda$, $R$, $p$ and $b$ result in a geometrically
ergodic sampler. One should choose the first three as discussed above, and
$b$ should be chosen to be small, say 0.1 or 0.01.

% All of our results were proven in $\mathbb{R}^k$, only assuming that
%$k$ is
% a positive integer. Hence the results hold in all dimensions.
% CHARLIE: not here, way to late to be helpful (if needed at all)

Admittedly, our methods do not guarantee geometric ergodicity
without any theoretical analysis. Users must understand
the tail behavior of the target distribution
in order to select the correct transformation. For
distributions with well behaved tails, this analysis may be easy, as in
our examples. We can say that our methods are no more difficult to apply
than the current state of the art [\citet{jarner-hansen}] and are applicable
to a much larger class of models.

\begin{appendix}\label{app}

%s5 #&#
\section{Isomorphic Markov chains} \label{secisomorphic}

We say measurable spaces are \emph{isomorphic} if there is an invertible
bimeasurable mapping between them ($h$ bimeasurable means
both $h$ and $h^{- 1}$ are measurable). We say probability spaces
$(S, \mathcal{A}, P)$ and $(T, \mathcal{B}, Q)$ are \emph{isomorphic}
if there is an invertible bimeasurable mapping $h \dvtx S \to T$ such that
$P = Q \circ h$, meaning
\[
P(A) = Q \bigl(h(A) \bigr),\qquad  A \in\mathcal{A},
\]
which also implies $Q = P \circ h^{- 1}$. We say Markov chains on state
spaces $(S, \mathcal{A})$ and $(T, \mathcal{B})$ are \emph
{isomorphic} if
there is an invertible bimeasurable mapping $h \dvtx S \to T$ such that the
corresponding initial distributions $\mu$ and $\nu$ and the transition
probability kernels $P$ and $Q$ satisfy $\mu= \nu\circ h$ and
%
%
%e30 #&#
\begin{equation}
\label{eqpq} P(x, A) = Q \bigl( h(x), h(A) \bigr),\qquad  x \in S \opand A \in
\mathcal{A}.
\end{equation}
By the change-of-variable theorem for measures, \eqref{eqpq} implies
%
%
%e31 #&#
\begin{equation}
\label{eqpq-too} P^n(x, A) = Q^n \bigl( h(x), h(A)
\bigr),\qquad  n \in\nats\opand x \in S \opand A \in\mathcal{A}.
\end{equation}
It follows that $P$ has an irreducibility measure if and only if $Q$
has an irreducibility measure.
It also follows from the change-of-variable theorem that
$\eta$ is an invariant measure for $P$ if and only if $\eta\circ h^{- 1}$
is an invariant measure for $Q$. Thus $P$ is null recurrent if and only
if $Q$ is, and $P$ is positive recurrent if and only if $Q$ is.
Also $P$ is reversible with respect to $\eta$ if and only if $Q$ is
reversible with respect to $\eta\circ h^{- 1}$.

For Harris recurrence we use the criterion that a recurrent Markov chain
is Harris if and only if every bounded harmonic function is constant
[\citet{nummelin}, Theorem~3.8 combined with his Proposition~3.9 and
Theorem~8.0.1 of \citet{meyn-tweedie}]. A function $g$ is \emph{harmonic}
for a kernel $P$ if $g = P g$, meaning
\[
g(x) = \int P(x, d y) g(y),\qquad  x \in S.
\]
It is clear that $g$ is harmonic for $P$ if and only if $g \circ h^{-
1}$ is
harmonic for $Q$. Thus~$P$ is Harris recurrent if and only
if $Q$ is.

Suppose $P$ is irreducible and periodic.
This means [\citet{meyn-tweedie}, Proposition~5.4.1] there are
disjoint sets
$D_0$, $\ldots,$ $D_{d - 1}$ with $d \ge2$ that are a partition of $S$
such that
\[
P(x, D_{{i + 1 \ \mathrm{mod}\  d}}) = 1,\qquad  x \in D_i, \ i = 0, \ldots, d - 1.
\]
But then
\[
Q \bigl(y, h^{- 1}(D_{{i + 1 \ \mathrm{mod}\  d}}) \bigr) = 1,\qquad  y \in
h^{- 1}(D_i), \ i = 0, \ldots, d - 1,
\]
and the sets $h^{- 1}(D_i)$ partition $T$, so $Q$ is also periodic.
Thus isomorphic irreducible Markov chains are both periodic or both aperiodic.

Finally suppose $\pi$ is an invariant probability measure for $P$,
and $\mu$ is any probability measure on the state space.
Then $\psi= \pi\circ h^{- 1}$ is an invariant probability measure for $Q$,
and it is clear that
\[
\bigl\Vert\pi- \mu P^{n} \bigl\Vert = \bigl\Vert\psi- \nu Q^n \bigr\Vert,\qquad n
\in \nats,
\]
where $\Vert\fatdot \Vert$ denotes total variation norm
and $\nu= \mu\circ h^{- 1}$. A Markov chain is \emph{geometrically
ergodic} if there exists a nonnegative-real-valued function $M$
and constant $r<1$
such that
%
%
%e32 #&#
\begin{equation}
\label{eqgeometrically-ergodic-definition}\bigl \Vert P^n(x, \cdot) - \pi
(\cdot) \bigr\Vert \leq M(x) r^n,\qquad \mbox{for all $x$}
\end{equation}
[\citet{meyn-tweedie}, Chapter~15].
If $M$ is bounded, then the Markov chain is \emph{uniformly ergodic}
[\citet{meyn-tweedie}, Chapter~16].
If \eqref{eqgeometrically-ergodic-definition} holds with $r^n$
replaced by
$n^r$ for some $r < 0$, then the Markov chain is \emph{polynomially ergodic}
[\citet{jarner-roberts-polynomial}].
Thus, if a Markov chain is polynomially ergodic, geometrically ergodic, or
uniformly ergodic,
then any isomorphic Markov chain has the same property.

The following summarizes the discussion in this appendix.
%
%
%th6 #&#
\begin{thmm}[(Isomorphic Markov chains)]
If a Markov chain has one of the following properties,
irreducibility,
reversibility,
null recurrence,
positive recurrence,
Harris recurrence,
aperiodicity,
polynomial ergodicity,
geometric ergodicity,
uniform ergodicity,
then so does any isomorphic Markov chain.
\end{thmm}

%s6 #&#
\section{\texorpdfstring{Proof of Lemma \lowercase{\protect\ref{lemf-conds}}}{Proof of Lemma 1}}
\label{secproof-lemma-isotropic}
% \lowercase necessary because the document class is uppercasing section
% titles and this uppercases the ref ??? WTF TeX ???

That $f$ is a diffeomorphism follows from the inverse function theorem
\[
\frac{d f^{- 1}(t)}{d t} = \frac{1}{f'(s)} \qquad\mbox{whenever $t = f(s)$}
\]
and \eqref{condf-first}. It is clear from \eqref{eqh} that
$\vert h(\gamma) \vert = f(\vert\gamma \vert)$ for all $\gamma$,
from which
\eqref{eqgamma-ratio}, \eqref{eqh-inverse} and the invertibility of
$h$ follow.

Now for $\gamma\neq0$ we have
\[
\frac{\partial}{\partial\gamma_k} \Biggl(\sum_{i = 1}^d
\gamma_i^2 \Biggr)^{1 / 2} = \Biggl(\sum
_{i = 1}^d \gamma_i^2
\Biggr)^{- 1 / 2} \gamma_k
\]
so
\[
\nabla\vert\gamma \vert = \frac{\gamma^T}{\vert\gamma \vert},
\]
and now \eqref{eqdee-h-nonzero} follows straightforwardly from
\eqref{eqh}, and it is clear that $h$ is continuously
differentiable everywhere except perhaps at zero and similarly for
$h^{- 1}$.

The term in square brackets on the right-hand side of
\eqref{eqdee-h-nonzero} goes to zero as $\vert\gamma \vert \to0$
by the
definition of derivative and that the term that multiplies it is bounded,
thus, if we can show \eqref{eqdee-h-zero}, then $\nabla h$ is also
continuous at zero. By the definition of derivative, what must be shown
to prove \eqref{eqdee-h-zero} is that
\begin{eqnarray*}
\frac{h(\gamma) - f'(0) \gamma}{\vert\gamma \vert} & =& \frac
{f(\vert\gamma \vert)
({\gamma}/{\vert\gamma \vert})
- f'(0) \gamma}{\vert\gamma \vert}
\\
& =& \biggl[ \frac{f(\vert\gamma \vert)}{\vert\gamma \vert} - f'(0) \biggr] \frac{\gamma}{\vert\gamma \vert}
\end{eqnarray*}
converges to zero as $\gamma\to0$. Since the term in square brackets
converges to zero by the definition of derivative and $\gamma/
\vert\gamma \vert$ is bounded, this proves \eqref{eqdee-h-zero}.
Since the
formulas for $h$ and $h^{- 1}$ have the same form, this shows $h$ is a
diffeomorphism.

The determinant of a symmetric matrix is the product of its eigenvalues
[\citet{harville}, Theorem~21.6.1]. First, $\gamma$ is an eigenvector of
$\nabla h(\gamma)$ with eigenvalue $f'(\vert\gamma \vert)$. Second, any
vector $v$ orthogonal to $\gamma$ is also an eigenvector of $\nabla
h(\gamma)$ with eigenvalue $f(\vert\gamma \vert) / \vert\gamma
\vert$ when
$\gamma\neq0$ and eigenvalue $f'(0)$ when $\gamma= 0$. Since the
subspace orthogonal to $\gamma$ has dimension $k - 1$, the multiplicity
of the second kind of eigenvalue is $k - 1$. This proves
\eqref{eqdet-dee-h}.

For $\gamma\neq0$ we have
%
%
%e33 #&#
%e34 #&#
\begin{eqnarray}
\label{eqdee-det-h}\qquad \nabla\det \bigl( \nabla h(\gamma) \bigr) &=&
f^{\prime\prime}\bigl(\vert\gamma \vert\bigr) \biggl( \frac{f(\vert
\gamma \vert)}{\vert\gamma \vert}
\biggr)^{k-1} \frac{\gamma^T}{\vert\gamma \vert}
\nonumber
\\[-8pt]
\\[-8pt]
\nonumber
&&{}+ (k-1) f^\prime\bigl(\vert\gamma \vert\bigr) \biggl( \frac{f(\vert
\gamma \vert)}{\vert\gamma \vert}
\biggr)^{k-2} \biggl[ \frac{f'(\vert\gamma \vert)}{\vert\gamma
\vert} - \frac{f(\vert\gamma \vert)}{\vert\gamma \vert^2} \biggr]
\frac{\gamma^T}{\vert\gamma \vert}.
\end{eqnarray}
Since \eqref{eqdet-dee-h} depends on $\gamma$ only through
$\vert\gamma \vert$, it has circular contours, and we must have
%
%
%e35 #&#
\begin{equation}
\label{eqdee-det-h-zero} \nabla\det \bigl( \nabla h(0) \bigr) = 0
\end{equation}
if the derivative exists. We claim the derivative
\eqref{eqdee-det-h-zero} does exist, and \eqref{eqdet-dee-h} is
continuously differentiable under the ``additional assumptions'' about
second derivatives of $f$ of the lemma. To prove this claim we need to
first show that \eqref{eqdee-det-h} converges to zero as $\gamma\to0$
and second show that \eqref{eqdee-det-h-zero} is the derivative at zero.

Except for the behavior of the term in square brackets, the limit of
\eqref{eqdee-det-h} is obvious from $f(s) / s \to f'(0)$ as $s \to0$
and $\gamma/ \vert\gamma \vert$ being bounded. For the term in square
brackets we use Taylor's theorem [\citet{stromberg}, Theorem~4.34]
\begin{eqnarray*}
f(s) & = &c s + o \bigl(s^2 \bigr),
\\
f'(s) & =& c + o(s),
\end{eqnarray*}
where $c = f'(0)$, so
\[
\frac{f'(s)}{s} - \frac{f(s)}{s^2} = o(1),
\]
and the term in square brackets in \eqref{eqdee-det-h} goes to zero as
$\gamma\to0$ proving that all of~\eqref{eqdee-det-h} goes to zero as
$\gamma\to0$.

What must be shown to establish \eqref{eqdee-det-h-zero} is that
\begin{eqnarray*}
\frac{\det(\nabla h(\gamma)) - \det(\nabla h(0))}{\vert\gamma
\vert}  = \frac{1}{\vert\gamma \vert} \biggl[ f'\bigl(\vert\gamma
\vert\bigr) \biggl[ \frac{f(\vert\gamma \vert)}{\vert\gamma \vert} \biggr]^{k - 1} - \bigl[
f'(0) \bigr]^k \biggr]
\end{eqnarray*}
converges to zero as $\gamma\to0$. Applying L'Hospital's rule, we have
\begin{eqnarray*}
&&\lim_{s \downarrow0} \frac{f'(s) [ {f(s)}/{s} ]^{k - 1}
- [ f'(0) ]^k
}{s}
\\
&&\qquad= \lim_{s \downarrow0} \biggl[ f''(s) \biggl[
\frac{f(s)}{s} \biggr]^{k - 1} + f'(s) (k - 1) \biggl[
\frac{f(s)}{s} \biggr]^{k - 2} \biggl( \frac{f'(s)}{s} -
\frac{f(s)}{s^2} \biggr) \biggr],
\end{eqnarray*}
and we have already shown that the limit on the right-hand side is
zero.

%s7 #&#
\section{\texorpdfstring{Proofs from Section \lowercase{\protect\ref{seclighter-tails}}}{Proofs from Section 2.3}}
\label{secproof-induced}

Before we prove Theorem~\ref{thmexponential-to-super} we need two
additional lemmas.

%
%le2 #&#
\begin{lem}\label{lemh-leftover}
Let $h$ be defined by \eqref{eqh} and \eqref{eqf1}. Then
%
%
%e36 #&#
\begin{equation}
\label{eqh-leftover} \lim_{\vert\gamma \vert \rightarrow\infty} \frac
{\gamma}{\vert\gamma \vert} \cdot\nabla\log\det
\bigl( \nabla h(\gamma) \bigr) = 0,
\end{equation}
where the dot indicates inner product.
\end{lem}
\begin{pf}
Recalling the value of $\det( \nabla h(\gamma) )$
for $\gamma\ne0$ from
\eqref{eqdet-dee-h} we can rewrite the dot product in
\eqref{eqh-leftover} as
%
%
%e37 #&#
\begin{equation}
\label{eqgamma-dot-deth} \frac{f^{\prime\prime}(\vert\gamma \vert
)}{f^\prime(\vert\gamma \vert)} + (k - 1 ) \biggl(
\frac{f^\prime
(\vert\gamma \vert) }{f(\vert\gamma \vert)} - \frac{1}{\vert\gamma \vert} \biggr).
\end{equation}
From \eqref{eqf1} for $\vert\gamma \vert > R$ we have
%
%
%e38 #&#
%e39 #&#
\begin{eqnarray}
f'(x) &= 1 + p(x - R)^{p-1} ,\label{eqf1-prime}
\\
f''(x) &= p(p-1) (x - R)^{p-2}
\label{eqf1-double-prime}
\end{eqnarray}
and, plugging these into \eqref{eqgamma-dot-deth}, we see that,
because $p > 2$, all terms in \eqref{eqgamma-dot-deth}
go to zero like $\vert\gamma \vert^{- 1}$ as $\vert\gamma \vert
\to\infty$.
\end{pf}

%
%le3 #&#
\begin{lem}\label{lemf-results}
Under the assumptions of Lemma~\ref{lemf-conds},
%
%
%e40 #&#
\begin{eqnarray}
\label{eqdee-h-gamma} \nabla h(\gamma) \gamma&=& f'\bigl(\vert\gamma
\vert\bigr) \gamma, \qquad\gamma\in\real^k,
\\
\label{eqnh-nh} \bigl[ \nabla h(\gamma) \bigr]^2 &=&
\frac{f(\vert\gamma \vert)^2}{\vert\gamma \vert^2} \I_k + \biggl[ f^\prime\bigl(\vert\gamma
\vert\bigr)^2 - \frac{f(\vert\gamma \vert
)^2}{\vert\gamma \vert^2} \biggr] \frac{\gamma\gamma^T}{\vert\gamma \vert^2},\qquad \gamma
\neq0,
\end{eqnarray}
$\nabla h(\gamma)$ being a symmetric matrix, and
%
%
%e42 #&#
%e43 #&#
\begin{eqnarray}
\label{eqx-nh-nh-x}  x^T \bigl[ \nabla h(\gamma) \bigr]^2
x
 = \frac{f(\vert\gamma \vert)^2}{\vert\gamma \vert^2} \vert x \vert^2 + \biggl[ f^\prime\bigl(
\vert\gamma \vert\bigr)^2 - \frac{f(\vert\gamma \vert
)^2}{\vert\gamma \vert^2} \biggr] \biggl(
\frac{h(\gamma) \cdot x}{\vert h(\gamma) \vert} \biggr)^2,
\nonumber
\\[-8pt]
\\[-8pt]
\eqntext{x \in\real^k, \gamma\neq0.}
\end{eqnarray}
\end{lem}

\begin{pf}
From \eqref{eqdee-h-nonzero} and \eqref{eqdee-h-zero}, we
straightforwardly obtain \eqref{eqdee-h-gamma} and for $\gamma\neq0$
%
%
%e44 #&#
%e45 #&#
\begin{eqnarray}
\label{eqnh-nh1} \bigl[ \nabla h(\gamma) \bigr]^2 &=& \nabla h(
\gamma) \biggl( \frac{f(\vert\gamma \vert)}{\vert\gamma \vert} \I_k + \biggl[ \frac{f^\prime(\vert\gamma \vert)}{\vert\gamma \vert^2} -
\frac{f(\vert\gamma \vert)}{\vert\gamma \vert^3} \biggr] \gamma \gamma^T \biggr)
\nonumber
\\[-8pt]
\\[-8pt]
\nonumber
&= &\frac{f(\vert\gamma \vert)}{\vert\gamma \vert} \nabla h(\gamma) + \biggl[ \frac{f^\prime(\vert\gamma \vert)^2}{\vert\gamma \vert^2} -
\frac{f(\vert\gamma \vert)f^\prime(\vert\gamma \vert)}{\vert
\gamma \vert^3} \biggr] \gamma\gamma^T
\end{eqnarray}
and
\[
\frac{f(\vert\gamma \vert)}{\vert\gamma \vert} \nabla h(\gamma) = \frac
{f(\vert\gamma \vert)^2}{\vert\gamma \vert^2} \I_k +
\biggl[ \frac{f^\prime(\vert\gamma \vert)f(\vert\gamma \vert
)}{\vert\gamma \vert^3} - \frac{f(\vert\gamma \vert)^2}{\vert
\gamma \vert^4} \biggr] \gamma
\gamma^T,
\]
which plugged into \eqref{eqnh-nh1} gives \eqref{eqnh-nh}, and
\eqref{eqx-nh-nh-x} is straightforward from \eqref{eqnh-nh}.
\end{pf}

\begin{pf*}{Proof of Theorem~\ref{thmexponential-to-super}}
Since $\nabla h(\gamma)$ is a symmetric matrix, it follows
from~\eqref{eqdensity-trans-dee-log} that
\[
\gamma\fatdot\nabla\log\pi_\gamma(\gamma) = \nabla h(\gamma) \gamma
\fatdot\log\pi_\beta \bigl( h(\gamma) \bigr) + \gamma\fatdot \nabla\log
\det \bigl( \nabla h(\gamma) \bigr).
\]
Hence we can bound \eqref{eqtails}
by the sum of
%
%
%e46 #&#
\begin{equation}
\label{eqexp-limsup-main} \limsup_{\vert\gamma \vert \rightarrow
\infty
} \frac{\nabla h(\gamma) \gamma}{\vert\gamma \vert} \fatdot
\nabla\log\pi_\beta \bigl( h(\gamma) \bigr)
\end{equation}
and
%
%
%e47 #&#
\begin{equation}
\label{eqexp-limsup-leftover} \limsup_{\vert\gamma \vert
\rightarrow
\infty} \frac{\gamma}{\vert\gamma \vert} \fatdot
\nabla\log\det \bigl( \nabla h(\gamma) \bigr).
\end{equation}
It follows from \eqref{eqgamma-ratio} and \eqref{eqdee-h-gamma}
that for large $\vert\gamma \vert$ the dot product in
\eqref{eqexp-limsup-main} can be rewritten as
%
%
%e48 #&#
\begin{equation}
\label{eqfred} f'\bigl(\vert\gamma \vert\bigr) \frac{h(\gamma)}{\vert
h(\gamma ) \vert} \fatdot
\nabla \log\pi_\beta \bigl( h(\gamma) \bigr).
\end{equation}
Since $f'(\vert\gamma \vert)$ is always positive, and
$\pi_\beta$ is exponentially light, there is an $\varepsilon> 0$
such that \eqref{eqfred} is bounded above by
$- f_1'(\vert\gamma \vert) \varepsilon$.
It is clear that $f'(\vert\gamma \vert) \to\infty$ as $\vert
\gamma \vert
\to
\infty$,
so \eqref{eqexp-limsup-main} is equal to $-\infty$. It
follows from Lemma~\ref{lemh-leftover} that
\eqref{eqexp-limsup-leftover} is equal to zero, so
\eqref{eqtails} is equal to $-\infty$ and $\pi_\gamma$ is a
super-exponentially light density.
\end{pf*}

Before we prove Theorem~\ref{thmsub-exponential-to-exponential} we
need a
lemma.

%
%le4 #&#
\begin{lem}\label{lemdet-f2}
Let $h$ be defined by \eqref{eqh} and \eqref{eqf2}. Then
%
%
%e49 #&#
\begin{equation}
\label{eqdet-f2} \limsup_{\vert\gamma \vert \rightarrow\infty} \frac
{\gamma}{\vert\gamma \vert} \cdot\nabla\log\det
\bigl( \nabla h(\gamma) \bigr) = bk,
\end{equation}
where the dot indicates inner product.
\end{lem}
\begin{pf}
As in in the proof of Lemma~\ref{lemh-leftover}, the dot product in
\eqref{eqdet-f2}
can be written as \eqref{eqgamma-dot-deth}.
Clearly, $(k-1) /
\vert\gamma \vert$ goes to zero as $\vert\gamma \vert$ goes to
infinity. Hence,
\eqref{eqdet-f2} is equal to
%
%
%e50 #&#
\begin{equation}
\label{eqsally} \limsup_{x \rightarrow\infty} \biggl[ \frac
{f^{\prime\prime}(x)}{f^\prime(x)} + (k-1)
\frac{f^\prime(x)}{f(x)} \biggr]
\end{equation}
if the limit exists. For $x > 1/b$, it follows from
\eqref{eqf2} that
\begin{eqnarray*}
f'(x) &=& be^{bx},
\\
f''(x) &=& b^2e^{bx}
\end{eqnarray*}
and plugging these into \eqref{eqsally} gives
\[
\limsup_{x \rightarrow\infty} \biggl[ \frac{b^2e^{bx}}{be^{bx}} + (k-1) \frac{be^{bx}}{e^{bx} - e/3}
\biggr],
\]
which equals $bk$.
\end{pf}

\begin{pf*}{Proof of Theorem~\ref{thmsub-exponential-to-exponential}}
As in the proof of Theorem~\ref{thmexponential-to-super}, \eqref{eqtails}
can be rewritten as the sum
of \eqref{eqexp-limsup-main} and \eqref{eqexp-limsup-leftover},
and for large $\vert\gamma \vert$ the dot product in
\eqref{eqexp-limsup-main} can be rewritten as \eqref{eqfred}.
By \eqref{condsub-bounded}
and the fact that $\vert h(\gamma) \vert = f(\vert\gamma \vert)$,
\eqref{eqfred} is bounded above
\[
\limsup_{\vert\gamma \vert \rightarrow\infty} \biggl( - \alpha \frac
{f'(\vert\gamma \vert)}{f(\vert\gamma \vert)} \biggr),
\]
which when $f$ is given by \eqref{eqf2} is equal to $-b \alpha$.
It follows that the limit superior in~\eqref{eqtails} is bounded above by
$-b(\alpha- k)$. Since $\alpha> k$, this upper bound is less than~0,
so $\pi_\gamma$ is exponentially light.
\end{pf*}

%s8 #&#
\section{\texorpdfstring{Proofs from Section \lowercase{\protect\ref{seccurvature}}}{Proofs from Section 2.4}}
\label{secproof-curvature}

Some lemmas are needed to prove the curvature conditions for exponentially
light densities.

%
%le5 #&#
\begin{lem}\label{lemlimsup-magnitude}
Let $\pi_\beta$ be an exponentially light density on $\real^k$, and let
$h$ be defined by \eqref{eqh} and \eqref{eqf1}. Then
%
%
%e51 #&#
\begin{equation}
\label{eqlim-pibeta-h} \bigl\vert\nabla\log\pi_\beta \bigl(h (\gamma) \bigr)
\nabla h(\gamma) \bigr\vert \rightarrow\infty\qquad \mbox{as } \vert \gamma \vert
\rightarrow \infty,
\end{equation}
%
% Charlie added \nabla in this equation to match Lemma 4 in the thesis,
% which this lemma seems to be
and $\pi_\gamma$ defined by \eqref{eqchange} has the property
%
%
%e52 #&#
\begin{equation}
\label{eqlimsup-magnitude} \lim_{\vert\gamma \vert \rightarrow
\infty} \frac{\vert\nabla\log\pi_\gamma(\gamma) \vert} {
\vert\nabla\log\pi_\beta( h(\gamma) ) \nabla h (\gamma) \vert} = 1.
\end{equation}
\end{lem}
\begin{pf}
The square of the left-hand side of \eqref{eqlim-pibeta-h} is,
by \eqref{eqx-nh-nh-x},
%
%
%e53 #&#
\begin{eqnarray}
\label{eqtodo} &&\frac{f(\vert\gamma \vert)^2}{\vert\gamma \vert^2} \bigl\vert\nabla\log\pi_\beta \bigl(h(
\gamma) \bigr) \bigr\vert^2
\nonumber
\\[-8pt]
\\[-8pt]
\nonumber
&&\qquad{}+ \biggl[ f^\prime\bigl(\vert\gamma
\vert\bigr)^2 - \frac{f(\vert\gamma
\vert)^2}{\vert\gamma \vert^2} \biggr] \biggl( \frac{h(\gamma) \cdot
\nabla\log\pi_\beta(h(\gamma) )
}{\vert h(\gamma) \vert}
\biggr)^2,
\end{eqnarray}
hence \eqref{eqlim-pibeta-h} holds if and only if \eqref{eqtodo} goes
to infinity.
Since the left-hand term of~\eqref{eqtodo} is nonnegative, it is
sufficient to show that the right-hand term goes to infinity to show that
all of \eqref{eqtodo} goes to infinity. By assumption $\pi_\beta$ is
exponentially light, and since $\vert h(\gamma) \vert = f(\vert
\gamma \vert)$,
there exists an $\varepsilon> 0$ and $M < \infty$ such that
\[
\frac{h(\gamma) \cdot
\nabla\log\pi_\beta(h(\gamma) )
}{\vert h(\gamma) \vert} \le- \varepsilon, \qquad \vert\gamma \vert \ge M.
\]
Thus in order to prove \eqref{eqtodo} goes to infinity as
$\vert\gamma \vert$ goes to infinity, it is sufficient to prove that
the term
in square brackets in \eqref{eqtodo} goes to infinity. Plugging in the
definitions of $f$ and $f'$ from \eqref{eqf1} and \eqref{eqf1-prime}
for large $x$, we obtain
\begin{eqnarray*}
f^\prime(x)^2 - \frac{f(x)^2}{x^2} & =& \bigl[1 + p (x -
R)^{p - 1} \bigr]^2 - \frac{[x + (x - R)^p]^2}{x^2}
\\
& =& \bigl(p^2 - 1 \bigr) x^{2 p - 2} + o \bigl(x^{2 p - 2}
\bigr),
\end{eqnarray*}
and since $p > 2$ by assumption, this goes to infinity as $x$ goes to
infinity; hence~\eqref{eqtodo} goes to infinity as $\vert\gamma
\vert$ goes
to infinity and \eqref{eqlim-pibeta-h} holds.

By \eqref{eqdensity-trans-dee-log}, showing that
\eqref{eqlimsup-magnitude} is true only requires showing that
%
%
%e54 #&#
\begin{equation}
\label{eqalice} \lim_{\vert\gamma \vert \rightarrow\infty} \frac{
\vert\nabla\log\det(\nabla h(\gamma) ) \vert } {
\vert\nabla\log\pi_\beta(h(\gamma) ) \nabla h(\gamma) \vert} = 0.
\end{equation}
It follows from \eqref{eqdet-dee-h} that for $\gamma\neq0$,
\[
\log\det \bigl(\nabla h(\gamma) \bigr) = \log f'\bigl(\vert\gamma
\vert\bigr) + (k-1) \log \biggl( \frac{f(\vert\gamma \vert)}{\vert\gamma \vert
} \biggr)
\]
and
%
%
%e55 #&#
\begin{equation}
\label{eqherman} \nabla\log\det \bigl(\nabla h(\gamma) \bigr) = \biggl(
\frac{f''(\vert\gamma \vert)}{f'(\vert\gamma \vert)} + (k-1) \biggl[ \frac{f'(\vert\gamma \vert)}{f(\vert\gamma \vert)} - \frac
{1}{\vert\gamma \vert} \biggr]
\biggr) \frac{\gamma^T}{\vert\gamma \vert}.
\end{equation}
Plugging in the definitions of $f$, $f'$ and $f''$ from \eqref{eqf1},
\eqref{eqf1-prime} and \eqref{eqf1-double-prime} for large~$x$, we see
that $f''(x) / f'(x)$ and $f'(x) / f(x)$ go to zero as $x$ goes to
infinity, and hence~\eqref{eqherman} goes to zero as $\vert\gamma
\vert$
goes to
infinity.
Hence the numerator in \eqref{eqalice} goes to zero.
By \eqref{eqlim-pibeta-h} the denominator in \eqref{eqalice} goes
to infinity, and hence \eqref{eqalice} holds.
\end{pf}

%
%le6 #&#
\begin{lem}\label{lemlimsup-ratio}
Let $\pi_\beta$ be an exponentially light density on $\real^k$, and
let $h$ be defined by \eqref{eqh} and \eqref{eqf1}. Then
$\pi_\gamma$ defined by \eqref{eqchange} has the property that
%
%
%e56 #&#
\begin{equation}
\label{eqlimsup-ratio-orig} \limsup_{\vert\gamma \vert \rightarrow
\infty} \frac{\gamma}{\vert\gamma \vert} \cdot
\frac{\nabla\pi_\gamma(\gamma)}{\vert\nabla\pi_\gamma(\gamma)
\vert}
\end{equation}
(which is the limit superior in the curvature condition) is bounded
above by
%
%
%e57 #&#
\begin{equation}
\label{eqlimsup-ratio-new} \limsup_{\vert\gamma \vert \rightarrow
\infty} f'\bigl(\vert\gamma
\vert\bigr) \frac{\gamma}{\vert\gamma \vert} \cdot\frac{\nabla\log\pi_\beta
(h(\gamma) )} {
\vert\nabla\log\pi_\beta(h(\gamma) ) \nabla h(\gamma) \vert},
\end{equation}
%
% Charlie: no \nabla in argument of pi_beta. Makes no sense as was.
%Right?
% Note that Lemma 6 of thesis also has this error.
where the dots in both equations denote inner products.
\end{lem}
\begin{pf}
We always assume that $\pi_\beta$ and $\pi_\gamma$ are positive
(Section~\ref{secpos-c1}), so we may take logs, obtaining
\[
\frac{\nabla\log\pi_\gamma(\gamma)}{\vert\nabla\log\pi_\gamma (\gamma) \vert} = \frac{\nabla\pi_\gamma(\gamma)}{\vert
\nabla\pi_\gamma(\gamma) \vert}.
\]
Thus \eqref{eqlimsup-ratio-orig} can be rewritten as
\[
\limsup_{\vert\gamma \vert \rightarrow\infty} \frac{\gamma
}{\vert\gamma \vert} \cdot\frac{\nabla\log\pi_\gamma(\gamma)} {
\vert\nabla\log\pi_\beta(h(\gamma) ) \nabla h(\gamma) \vert}
\frac{\vert\nabla\log\pi_\beta(h(\gamma) ) \nabla h(\gamma)
\vert} {
\vert\nabla\log\pi_\gamma(\gamma) \vert},
\]
and then we can use Lemma~\ref{lemlimsup-magnitude} as
\[
\limsup_{\vert\gamma \vert \rightarrow\infty} \frac{\gamma
}{\vert\gamma \vert} \cdot\frac{\nabla\log\pi_\gamma(\gamma)} {
\vert\nabla\log\pi_\beta(h(\gamma) ) \nabla h(\gamma) \vert}.
\]
If we expand $\nabla\log\pi_\gamma(\gamma)$ using
\eqref{eqdensity-trans-dee-log}, this is bounded above by
the sum of
%
%
%e58 #&#
\begin{equation}
\label{eqlimsup-ratio-newtarget} \limsup_{\vert\gamma \vert
\rightarrow
\infty} \frac{\gamma}{\vert\gamma \vert} \cdot
\frac{\nabla\log\pi_\beta(h(\gamma) ) \nabla h(\gamma)} {
\vert\nabla\log\pi_\beta(h(\gamma) ) \nabla h(\gamma) \vert}
\end{equation}
and
%
%
%e59 #&#
\begin{equation}
\label{eqlimsup-ratio-newtarget-too} \limsup_{\vert\gamma \vert
\rightarrow\infty} \frac{\gamma}{\vert\gamma \vert}
\cdot\frac{\nabla\log\det(
\nabla h(\gamma) )}{\vert\nabla \log\pi_\beta(h(\gamma) ) \nabla
h(\gamma) \vert}.
\end{equation}
It follows from Lemmas~\ref{lemh-leftover} and~\ref{lemlimsup-magnitude}
that \eqref{eqlimsup-ratio-newtarget-too} is zero. Hence the $\limsup
$ in
\eqref{eqlimsup-ratio-orig} is bounded above by
\eqref{eqlimsup-ratio-newtarget}, which is equal to
\eqref{eqlimsup-ratio-new} since $\nabla h(\gamma)$ is symmetric and
$\nabla h(\gamma) \gamma= f^\prime(\vert\gamma \vert) \gamma$.
% on rereading this, I found the last bit very confusing until I
%remembered
% that \nabla\log\pi_\beta\bigl(h(\gamma)\bigr)
% is to be interpreted as a row vector (!) so
% \gamma\cdot\nabla\log\pi_\beta\bigl(h(\gamma)\bigr) \nabla h(
% does indeed have the interpretation
% \nabla\log\pi_\beta\bigl(h(\gamma)\bigr) \nabla h(\gamma)} \gamma
% confusing! Do you think we need some comment about this?
\end{pf}

%
%le7 #&#
\begin{lem}\label{lemweird-fraction}
Let $a(\gamma)$ and $b(\gamma)$ be functions such that both $a$ and $b$
are positive and bounded away from zero and infinity as $\vert\gamma
\vert$
goes to infinity. Then for $f$ from~\eqref{eqf1}, the fraction
%
%
%e60 #&#
\begin{equation}
\label{eqweird-fraction} f'\bigl(\vert\gamma \vert\bigr)^2 \Big/
\biggl( \frac{f(\vert\gamma \vert)^2}{\vert\gamma \vert^2} a(\gamma) + \biggl[ f'\bigl(\vert\gamma
\vert\bigr)^2 - \frac{f(\vert\gamma \vert)^2}{\vert
\gamma \vert^2} \biggr] b(\gamma) \biggr)
\end{equation}
is positive and bounded away from zero and infinity as $\vert\gamma
\vert$
goes to
infinity.
\end{lem}
\begin{pf}
The reciprocal of \eqref{eqweird-fraction} is
\[
\frac{f(\vert\gamma \vert)^2}{f'(\vert\gamma \vert)^2 \vert
\gamma \vert^2} a(\gamma) + \biggl[ 1 - \frac{f(\vert\gamma \vert)^2}{f'(\vert
\gamma \vert)^2\vert\gamma \vert^2} \biggr] b(\gamma).
\]
Since $a(\gamma)$ and $b(\gamma)$ are both positive and bounded away from
zero and infinity for large $\vert\gamma \vert$, it is sufficient to show
that
%
%
%e61 #&#
\begin{equation}
\label{eqcurv-bdd-lim-simplified} \frac{f(x)^2}{f'(x)^2 x^2}
\end{equation}
is bounded away from zero and one for large $x$.
For large $x$, it follows from \eqref{eqf1} and \eqref{eqf1-prime} that
\eqref{eqcurv-bdd-lim-simplified} is equal to
\[
\frac{[x + (x-R)^p]^2}{[1+p(x-R)^{p-1}]^2 x^2},
\]
which converges to $1 / p^2$ as $x \to\infty$. Since we assume $p > 2$,
we are done.
\end{pf}

\begin{pf*}{Proof of Theorem~\ref{thmcurvature-exponential}}
First, assume that condition (i) holds.
By Lemma~\ref{lemlimsup-ratio}, it is enough to show that
\eqref{eqlimsup-ratio-new} is less than zero, and
\eqref{eqlimsup-ratio-new} is equal to, using
\eqref{eqgamma-ratio},
%
%
%e62 #&#
\begin{equation}
\label{eqrearrange} \limsup_{\vert\gamma \vert \rightarrow\infty} \frac{\vert\nabla\log\pi_\beta(h(\gamma) ) \vert
f'(\vert\gamma \vert)} {
\vert\nabla\log\pi_\beta(h(\gamma) ) \nabla h(\gamma) \vert} \frac
{h(\gamma)}{\vert h(\gamma) \vert}
\cdot\frac{\nabla\log\pi_\beta(h(\gamma) )} {
\vert\nabla\log\pi_\beta(h(\gamma) ) \vert}.
\end{equation}
Since $\pi_\beta$ satisfies condition \eqref{condjh-curvature},
there is
an $\varepsilon> 0$ such that \eqref{eqrearrange} is bounded above by
%
%
%e63 #&#
\begin{equation}
\label{eqdel-pi-curv-target} \limsup_{\vert\gamma \vert
\rightarrow
\infty} \frac{\vert\nabla\log\pi_\beta(h(\gamma) ) \vert
f^\prime(\vert\gamma \vert)} {
\vert\nabla\log\pi_\beta(h(\gamma) ) \nabla h(\gamma) \vert} (-
\varepsilon).
\end{equation}
Because $f'(\vert\gamma \vert)$ is strictly positive, the fraction in
\eqref{eqdel-pi-curv-target}
is strictly positive for large~$\vert\gamma \vert$, hence
showing that this fraction's square is bounded away from zero is
enough to show that \eqref{eqdel-pi-curv-target} is less than zero,
and condition \eqref{condjh-curvature} holds. Let
\[
a(\gamma) = \frac{\vert\nabla\log \pi_\beta(h(\gamma) ) \vert^2}{\vert\nabla\log \pi_\beta(h(\gamma) ) \vert^2} = 1
\]
and
\[
b(\gamma) = \biggl( \frac{\nabla\log\pi_\beta(h(\gamma) ) \cdot
h(\gamma)} {
\vert\nabla\log\pi_\beta(h(\gamma) ) \vert
\vert h(\gamma) \vert
} \biggr)^2.
\]
Then, using \eqref{eqx-nh-nh-x} as in deriving \eqref{eqtodo},
the square of the fraction in \eqref{eqdel-pi-curv-target} is equal to~\eqref{eqweird-fraction}.
The Cauchy--Schwarz inequality
bounds $b(\gamma)$ above by one, and condition~\eqref{condjh-curvature}
bounds $b(\gamma)$ away from zero. So by Lemma~\ref{lemweird-fraction}
the square of the fraction in~\eqref{eqdel-pi-curv-target} is positive
and bounded away from zero as $\vert\gamma \vert$ goes to infinity. Because
this fraction itself is positive, it must also be bounded away from zero
as $\vert\gamma \vert$ goes to infinity. Hence the $\limsup$ in
\eqref{eqdel-pi-curv-target} is negative and condition
\eqref{condjh-curvature} holds for $\pi_\gamma$.

Now assume that condition (ii) holds
and $\pi_\beta$ is exponentially light, that is, there exist a $\beta_0 > 0$,
$\varepsilon> 0$ and $M_1 > M_2 > 0$ such that for $\vert\beta \vert >
\beta_0$,
\[
\frac{\beta}{\vert\beta \vert} \cdot\nabla\log\pi_\beta(\beta ) < - \varepsilon
\]
and
\[
M_2 < \bigl\vert\nabla\log\pi_\beta(\beta) \bigr\vert <
M_1.
\]
It follows that $1 / \vert\nabla\log\pi_\beta(\beta) \vert > 1 /
M_1$ so
$\pi_\beta$ satisfies condition~(i).
\end{pf*}

\begin{pf*}{Proof of Theorem~\ref{thmcurvature-sub-exponential}}
By \eqref{eqdensity-trans-dee-log} and the triangle inequality,
$\vert\nabla\log\pi_\gamma(\gamma) \vert$ is bounded above by
the sum
%
%
%e64 #&#
\begin{equation}
\label{eqnorm-bdd-sum} \bigl\vert\nabla\log\pi_\beta \bigl( h(\gamma )
\bigr) \nabla h(\gamma) \bigr\vert + \bigl\vert\nabla\log\det \bigl( h(\gamma) \bigr) \bigr\vert.
\end{equation}
Hence it is sufficient to show that both of these terms are bounded as
$\vert\gamma \vert$ goes to infinity.

It follows from \eqref{eqherman} that the right-hand
term in \eqref{eqnorm-bdd-sum} is equal to
%
%
%e65 #&#
\begin{equation}
\label{eqherman-too} \frac{f^{\prime\prime}(\vert\gamma \vert
)}{f^\prime(\vert\gamma \vert)} + (k-1) \frac{f^\prime(\vert
\gamma \vert)}{f(\vert\gamma \vert)} - (k-1)
\frac{1}{\vert\gamma \vert}.
\end{equation}
For large $y$,
%
%
%e66 #&#
%e67 #&#
%e68 #&#
\begin{eqnarray}
f(y) & =& e^{by} - \frac{e}{3}, \label{eqf2-large-zero}
\\
f'(y) & =& be^{by}, \label{eqf2-large-one}
\\
f''(y) & = &b^2e^{by}
\label{eqf2-large-two}.
\end{eqnarray}
So \eqref{eqherman-too} is equal to
\[
b + b (k-1) \frac{e^{b \vert\gamma \vert}}{e^{b \vert\gamma \vert
} - e/3} - (k-1) \frac{1}{\vert\gamma \vert},
\]
which clearly converges to $b k$ as $\vert\gamma \vert$ goes to infinity,
so the
right-hand term in \eqref{eqnorm-bdd-sum} is bounded
for large $\vert\gamma \vert$.

It follows from \eqref{eqx-nh-nh-x} as in deriving \eqref{eqtodo}
and from \eqref{eqgamma-ratio} that the
square of the left-hand term in \eqref{eqnorm-bdd-sum} is equal to the
sum of
%
%
%e69 #&#
\begin{equation}
\label{eqbdd-tmp1} \frac{f(\vert\gamma \vert)^2}{\vert\gamma
\vert^2} \bigl\vert\nabla\log\pi_\beta \bigl(
h(\gamma) \bigr) \bigr\vert^2
\end{equation}
and
%
%
%e70 #&#
\begin{equation}
\label{eqbdd-tmp2} f'\bigl(\vert\gamma \vert\bigr)^2 \biggl[ 1 -
\frac{f(\vert\gamma \vert)^2}{\vert\gamma \vert^2f'(\vert
\gamma \vert)^2} \biggr] \biggl( \frac{h(\gamma)
\cdot
\nabla\log\pi_\beta(h(\gamma) )} {
\vert h(\gamma) \vert} \biggr)^2.
\end{equation}

It follows from \eqref{eqf2-large-zero} and \eqref{eqf2-large-one}
that the term in square brackets of \eqref{eqbdd-tmp2} is positive
and less than one for
large $\vert\gamma \vert$. Since the other two terms in \eqref{eqbdd-tmp2}
are squares, \eqref{eqbdd-tmp2}~is nonnegative for large
$\vert\gamma \vert$. Thus, applying the Cauchy--Schwarz inequality
to the
term in parentheses in \eqref{eqbdd-tmp2}, one bounds
\eqref{eqbdd-tmp2} above by
%
%
%e71 #&#
\begin{equation}
\label{eqbdd-tmp3} f'\bigl(\vert\gamma \vert\bigr)^2\bigl \vert\nabla
\log \pi_\beta \bigl( h(\gamma) \bigr) \bigr\vert^2.
\end{equation}

By $f(\vert\gamma \vert) = \vert h(\gamma) \vert$ and
by \eqref{condsub-magnitude-bdd}, for $\vert\gamma \vert$
large \eqref{eqbdd-tmp3} is bounded above by
\[
\alpha^2 \frac{f'(\vert\gamma \vert)^2}{f(\vert\gamma \vert)^2},
\]
which converges to $\alpha^2 b^2$ as $\vert\gamma \vert$ goes to infinity,
and that finishes the proof that \eqref{eqnorm-bdd-sum}
is bounded for large $\vert\gamma \vert$ and the proof of the theorem.
\end{pf*}
\end{appendix}

% imsref loaded by akundreckaite, 2012-11-15 12:15:50
% imsref loaded by akundreckaite, 2012-11-15 12:28:31
% imsref loaded by akundreckaite, 2012-11-15 12:46:21
%

%suskaldyti doi

\printaddresses


\begin{thebibliography}{48}
% BibTex style file: ims.bst, 2012-08-21
% Default style options (sort=0,type=number).
% Used options (sort=1,type=nameyear).

%b1 #&#
\bibitem[\protect\citeauthoryear{Barndorff-Nielsen}{1978}]{barndorff}
%
\begin{bbook}[mr]
\bauthor{\bsnm{Barndorff-Nielsen},~\bfnm{Ole}\binits{O.}}
(\byear{1978}).
\btitle{Information and Exponential Families in Statistical Theory}.
\bpublisher{Wiley}, \blocation{Chichester}.
\bid{mr={0489333}}
\bptok{imsref}%
\end{bbook}
%
\endbibitem

%b2 #&#
\bibitem[\protect\citeauthoryear{Brooks et~al.}{2011}]{brooks}
%
\begin{bbook}[mr]
\beditor{\bsnm{Brooks},~\bfnm{Steve}\binits{S.}},
\beditor{\bsnm{Gelman},~\bfnm{Andrew}\binits{A.}},
\beditor{\bsnm{Jones},~\bfnm{Galin~L.}\binits{G.~L.}} \AND
\beditor{\bsnm{Meng},~\bfnm{Xiao-Li}\binits{X.-L.}}, eds.
(\byear{2011}).
\btitle{Handbook of {M}arkov Chain {M}onte {C}arlo}.
\bseries{Chapman \& Hall/CRC Handbooks of Modern Statistical Methods}.
\bpublisher{CRC Press}, \blocation{Boca Raton, FL}.
\bid{doi={10.1201/b10905}, mr={2742422}}
\bptok{imsref}%
\end{bbook}
%
\endbibitem

%b3 #&#
\bibitem[\protect\citeauthoryear{Chan}{1993}]{chan}
%
\begin{barticle}[mr]
\bauthor{\bsnm{Chan},~\bfnm{K.~S.}\binits{K.~S.}}
(\byear{1993}).
\btitle{On the central limit theorem for an ergodic {M}arkov chain}.
\bjournal{Stochastic Process. Appl.}
\bvolume{47}
\bpages{113--117}.
\bid{doi={10.1016/0304-4149(93)90097-N}, issn={0304-4149}, mr={1232855}}
\bptok{imsref}%
\end{barticle}
%
\endbibitem

%b4 #&#
\bibitem[\protect\citeauthoryear{Chan and Geyer}{1994}]{chan-geyer}
%
\begin{barticle}[author]
\bauthor{\bsnm{Chan},~\bfnm{K.~S.}\binits{K.~S.}} \AND
\bauthor{\bsnm{Geyer},~\bfnm{C.~J.}\binits{C.~J.}}
(\byear{1994}).
\btitle{{C}omment on ``{M}arkov chains for exploring posterior
distributions.''}
\bjournal{Ann. Statist.}
\bvolume{22}
\bpages{1747--1758}.
\bptok{imsref}%
\end{barticle}
%
\endbibitem

%b5 #&#
\bibitem[\protect\citeauthoryear{Diaconis and
Ylvisaker}{1979}]{diaconis-ylvisaker}
%
\begin{barticle}[mr]
\bauthor{\bsnm{Diaconis},~\bfnm{Persi}\binits{P.}} \AND
\bauthor{\bsnm{Ylvisaker},~\bfnm{Donald}\binits{D.}}
(\byear{1979}).
\btitle{Conjugate priors for exponential families}.
\bjournal{Ann. Statist.}
\bvolume{7}
\bpages{269--281}.
\bid{issn={0090-5364}, mr={0520238}}
\bptok{imsref}%
\end{barticle}
%
\endbibitem

%b6 #&#
\bibitem[\protect\citeauthoryear{Flegal and Jones}{2010}]{flegal-jones}
%
\begin{barticle}[mr]
\bauthor{\bsnm{Flegal},~\bfnm{James~M.}\binits{J.~M.}} \AND
\bauthor{\bsnm{Jones},~\bfnm{Galin~L.}\binits{G.~L.}}
(\byear{2010}).
\btitle{Batch means and spectral variance estimators in {M}arkov chain {M}onte
{C}arlo}.
\bjournal{Ann. Statist.}
\bvolume{38}
\bpages{1034--1070}.
\bid{doi={10.1214/09-AOS735}, issn={0090-5364}, mr={2604704}}
\bptok{imsref}%
\end{barticle}
%
\endbibitem

%b7 #&#
\bibitem[\protect\citeauthoryear{Gelfand and Smith}{1990}]{gelfand-smith}
%
\begin{barticle}[mr]
\bauthor{\bsnm{Gelfand},~\bfnm{Alan~E.}\binits{A.~E.}} \AND
\bauthor{\bsnm{Smith},~\bfnm{Adrian F.~M.}\binits{A.~F.~M.}}
(\byear{1990}).
\btitle{Sampling-based approaches to calculating marginal densities}.
\bjournal{J. Amer. Statist. Assoc.}
\bvolume{85}
\bpages{398--409}.
\bid{issn={0162-1459}, mr={1141740}}
\bptok{imsref}%
\end{barticle}
%
\endbibitem

%b8 #&#
\bibitem[\protect\citeauthoryear{Geman and Geman}{1984}]{geman-geman}
%
\begin{barticle}[author]
\bauthor{\bsnm{Geman},~\bfnm{S.}\binits{S.}} \AND
\bauthor{\bsnm{Geman},~\bfnm{D.}\binits{D.}}
(\byear{1984}).
\btitle{{S}tochastic relaxtion, {G}ibbs distributions, and the {B}ayesian
restoration of images}.
\bjournal{IEEE Trans. Pattern Anal. Mach. Intell.}
\bvolume{6}
\bpages{721--741}.
\bptok{imsref}%
\end{barticle}
%
\endbibitem

%b9 #&#
\bibitem[\protect\citeauthoryear{Geyer}{1992}]{geyer-orange}
%
\begin{barticle}[author]
\bauthor{\bsnm{Geyer},~\bfnm{Charles~J.}\binits{C.~J.}}
(\byear{1992}).
\btitle{Practical {M}arkov chain {M}onte {C}arlo (with discussion)}.
\bjournal{Statist. Sci.}
\bvolume{7}
\bpages{473--511}.
\bptok{imsref}%
\end{barticle}
%
\endbibitem


%b11 #&#
\bibitem[\protect\citeauthoryear{Geyer}{2011}]{intro}
%
\begin{bincollection}[author]
\bauthor{\bsnm{Geyer},~\bfnm{Charles~J.}\binits{C.~J.}}
(\byear{2011}).
\btitle{Introduction to MCMC}.
In \bbooktitle{Handbook of {M}arkov Chain {M}onte {C}arlo}
(\beditor{\bfnm{S.~P.}\binits{S.~P.}~\bsnm{Brooks}},
\beditor{\bfnm{A.~E.}\binits{A.~E.}~\bsnm{Gelman}},
\beditor{\bfnm{G.~L.}\binits{G.~L.}~\bsnm{Jones}} \AND
\beditor{\bfnm{X.~L.}\binits{X.~L.}~\bsnm{Meng}}, eds.).
\bpublisher{Chapman \& Hall/CRC}, \blocation{Boca Raton}.
\bptok{imsref}%
\end{bincollection}
%
\endbibitem


%b10 #&#
\bibitem[\protect\citeauthoryear{Geyer and Johnson}{2012}]{mcmc-R-package}
%
\begin{bmisc}[author]
\bauthor{\bsnm{Geyer},~\bfnm{Charles~J.}\binits{C.~J.}}
\AND
\bauthor{\bsnm{Johnson},~\bfnm{L.~T.}\binits{L.~T.}}
(\byear{2012}).
\bhowpublished{mcmc: {M}arkov {C}hain {M}onte {C}arlo. R package
version~0.8. Available at \url{http://CRAN.R-project.org/package=mcmc}.}
\bptok{imsref}%
\end{bmisc}
%
\endbibitem

%b12 #&#
\bibitem[\protect\citeauthoryear{Geyer and M{\o}ller}{1994}]{geyer-moller}
%
\begin{barticle}[mr]
\bauthor{\bsnm{Geyer},~\bfnm{Charles~J.}\binits{C.~J.}} \AND
\bauthor{\bsnm{M{\o}ller},~\bfnm{Jesper}\binits{J.}}
(\byear{1994}).
\btitle{Simulation procedures and likelihood inference for spatial point
processes}.
\bjournal{Scand. J. Stat.}
\bvolume{21}
\bpages{359--373}.
\bid{issn={0303-6898}, mr={1310082}}
\bptok{imsref}%
\end{barticle}
%
\endbibitem

%b13 #&#
\bibitem[\protect\citeauthoryear{Gilks, Richardson and
Spiegelhalter}{1996}]{gilks}
%
\begin{bbook}[auto]
\beditor{\bsnm{Gilks},~\bfnm{W.~R.}\binits{W.~R.}},
\beditor{\bsnm{Richardson},~\bfnm{S.}\binits{S.}} \AND
\beditor{\bsnm{Spiegelhalter},~\bfnm{D.~J.}\binits{D.~J.}}, eds.
(\byear{1996}).
\btitle{Markov Chain {M}onte {C}arlo in Practice}.
\bseries{Interdisciplinary Statistics}.
\bpublisher{Chapman \& Hall}, \blocation{London}.
\bid{mr={1397966}}
\bptok{imsref}%
\end{bbook}
%
\endbibitem

%b14 #&#
\bibitem[\protect\citeauthoryear{Gordin and Lif{\v
{s}}ic}{1978}]{gordin-lifsic}
%
\begin{barticle}[mr]
\bauthor{\bsnm{Gordin},~\bfnm{M.~I.}\binits{M.~I.}} \AND
\bauthor{\bsnm{Lif{\v{s}}ic},~\bfnm{B.~A.}\binits{B.~A.}}
(\byear{1978}).
\btitle{Central limit theorem for stationary {M}arkov processes}.
\bjournal{Dokl. Akad. Nauk SSSR}
\bvolume{239}
\bpages{766--767}.
\bid{issn={0002-3264}, mr={0501277}}
\bptok{imsref}%
\end{barticle}
%
\endbibitem

%b15 #&#
\bibitem[\protect\citeauthoryear{Green}{1995}]{green}
%
\begin{barticle}[mr]
\bauthor{\bsnm{Green},~\bfnm{Peter~J.}\binits{P.~J.}}
(\byear{1995}).
\btitle{Reversible jump {M}arkov chain {M}onte {C}arlo computation and
{B}ayesian model determination}.
\bjournal{Biometrika}
\bvolume{82}
\bpages{711--732}.
\bid{doi={10.1093/biomet/82.4.711}, issn={0006-3444}, mr={1380810}}
\bptok{imsref}%
\end{barticle}
%
\endbibitem

%b16 #&#
\bibitem[\protect\citeauthoryear{Harville}{1997}]{harville}
%
\begin{bbook}[author]
\bauthor{\bsnm{Harville},~\bfnm{David~A.}\binits{D.~A.}}
(\byear{1997}).
\btitle{Matrix Algebra from a Statistician's Perspective}.
\bpublisher{Springer}, \blocation{New York}.
\bptok{imsref}%
\end{bbook}
%
\endbibitem

%b17 #&#
\bibitem[\protect\citeauthoryear{Hastings}{1970}]{hastings}
%
\begin{barticle}[author]
\bauthor{\bsnm{Hastings},~\bfnm{W.~K.}\binits{W.~K.}}
(\byear{1970}).
\btitle{{M}onte {C}arlo sampling methods using {M}arkov chains and their
applications}.
\bjournal{Biometrika}
\bvolume{57}
\bpages{97--109}.
\bptok{imsref}%
\end{barticle}
%
\endbibitem

%b18 #&#
\bibitem[\protect\citeauthoryear{Hobert and Geyer}{1998}]{hobert-geyer}
%
\begin{barticle}[mr]
\bauthor{\bsnm{Hobert},~\bfnm{James~P.}\binits{J.~P.}} \AND
\bauthor{\bsnm{Geyer},~\bfnm{Charles~J.}\binits{C.~J.}}
(\byear{1998}).
\btitle{Geometric ergodicity of {G}ibbs and block {G}ibbs samplers for a
hierarchical random effects model}.
\bjournal{J. Multivariate Anal.}
\bvolume{67}
\bpages{414--430}.
\bid{doi={10.1006/jmva.1998.1778}, issn={0047-259X}, mr={1659196}}
\bptok{imsref}%
\end{barticle}
%
\endbibitem

%b19 #&#
\bibitem[\protect\citeauthoryear{Jarner and Hansen}{2000}]{jarner-hansen}
%
\begin{barticle}[mr]
\bauthor{\bsnm{Jarner},~\bfnm{S{\o}ren~Fiig}\binits{S.~F.}} \AND
\bauthor{\bsnm{Hansen},~\bfnm{Ernst}\binits{E.}}
(\byear{2000}).
\btitle{Geometric ergodicity of {M}etropolis algorithms}.
\bjournal{Stochastic Process. Appl.}
\bvolume{85}
\bpages{341--361}.
\bid{doi={10.1016/S0304-4149(99)00082-4}, issn={0304-4149}, mr={1731030}}
\bptok{imsref}%
\end{barticle}
%
\endbibitem

%b20 #&#
\bibitem[\protect\citeauthoryear{Jarner and
Roberts}{2002}]{jarner-roberts-polynomial}
%
\begin{barticle}[mr]
\bauthor{\bsnm{Jarner},~\bfnm{S{\o}ren~F.}\binits{S.~F.}} \AND
\bauthor{\bsnm{Roberts},~\bfnm{Gareth~O.}\binits{G.~O.}}
(\byear{2002}).
\btitle{Polynomial convergence rates of {M}arkov chains}.
\bjournal{Ann. Appl. Probab.}
\bvolume{12}
\bpages{224--247}.
\bid{doi={10.1214/aoap/1015961162}, issn={1050-5164}, mr={1890063}}
\bptok{imsref}%
\end{barticle}
%
\endbibitem

%b21 #&#
\bibitem[\protect\citeauthoryear{Jarner and Roberts}{2007}]{jarner-roberts}
%
\begin{barticle}[mr]
\bauthor{\bsnm{Jarner},~\bfnm{S{\o}ren~F.}\binits{S.~F.}} \AND
\bauthor{\bsnm{Roberts},~\bfnm{Gareth~O.}\binits{G.~O.}}
(\byear{2007}).
\btitle{Convergence of heavy-tailed {M}onte {C}arlo {M}arkov chain algorithms}.
\bjournal{Scand. J. Stat.}
\bvolume{34}
\bpages{781--815}.
\bid{doi={10.1111/j.1467-9469.2007.00557.x}, issn={0303-6898}, mr={2396939}}
\bptok{imsref}%
\end{barticle}
%
\endbibitem

%b22 #&#
\bibitem[\protect\citeauthoryear{Jarner and Tweedie}{2003}]{jarner-tweedie}
%
\begin{barticle}[mr]
\bauthor{\bsnm{Jarner},~\bfnm{S{\o}ren~F.}\binits{S.~F.}} \AND
\bauthor{\bsnm{Tweedie},~\bfnm{Richard~L.}\binits{R.~L.}}
(\byear{2003}).
\btitle{Necessary conditions for geometric and polynomial ergodicity of
random-walk-type {M}arkov chains}.
\bjournal{Bernoulli}
\bvolume{9}
\bpages{559--578}.
\bid{doi={10.3150/bj/1066223269}, issn={1350-7265}, mr={1996270}}
\bptok{imsref}%
\end{barticle}
%
\endbibitem

%b23 #&#
\bibitem[\protect\citeauthoryear{Johnson and Jones}{2010}]{johnson-jones}
%
\begin{barticle}[mr]
\bauthor{\bsnm{Johnson},~\bfnm{Alicia~A.}\binits{A.~A.}} \AND
\bauthor{\bsnm{Jones},~\bfnm{Galin~L.}\binits{G.~L.}}
(\byear{2010}).
\btitle{Gibbs sampling for a {B}ayesian hierarchical general linear model}.
\bjournal{Electron. J. Stat.}
\bvolume{4}
\bpages{313--333}.
\bid{doi={10.1214/09-EJS515}, issn={1935-7524}, mr={2645487}}
\bptok{imsref}%
\end{barticle}
%
\endbibitem

%b24 #&#
\bibitem[\protect\citeauthoryear{Jones}{2004}]{jones}
%
\begin{barticle}[mr]
\bauthor{\bsnm{Jones},~\bfnm{Galin~L.}\binits{G.~L.}}
(\byear{2004}).
\btitle{On the {M}arkov chain central limit theorem}.
\bjournal{Probab. Surv.}
\bvolume{1}
\bpages{299--320}.
\bid{doi={10.1214/154957804100000051}, issn={1549-5787}, mr={2068475}}
\bptok{imsref}%
\end{barticle}
%
\endbibitem

%b25 #&#
\bibitem[\protect\citeauthoryear{Jones and Hobert}{2004}]{jones-hobert}
%
\begin{barticle}[mr]
\bauthor{\bsnm{Jones},~\bfnm{Galin~L.}\binits{G.~L.}} \AND
\bauthor{\bsnm{Hobert},~\bfnm{James~P.}\binits{J.~P.}}
(\byear{2004}).
\btitle{Sufficient burn-in for {G}ibbs samplers for a hierarchical random
effects model}.
\bjournal{Ann. Statist.}
\bvolume{32}
\bpages{784--817}.
\bid{doi={10.1214/009053604000000184}, issn={0090-5364}, mr={2060178}}
\bptok{imsref}%
\end{barticle}
%
\endbibitem

%b26 #&#
\bibitem[\protect\citeauthoryear{Kipnis and Varadhan}{1986}]{kipnis-varadhan}
%
\begin{barticle}[mr]
\bauthor{\bsnm{Kipnis},~\bfnm{C.}\binits{C.}} \AND
\bauthor{\bsnm{Varadhan},~\bfnm{S.~R.~S.}\binits{S.~R.~S.}}
(\byear{1986}).
\btitle{Central limit theorem for additive functionals of reversible {M}arkov
processes and applications to simple exclusions}.
\bjournal{Comm. Math. Phys.}
\bvolume{104}
\bpages{1--19}.
\bid{issn={0010-3616}, mr={0834478}}
\bptok{imsref}%
\end{barticle}
%
\endbibitem

%b27 #&#
\bibitem[\protect\citeauthoryear{{\L}atuszy{\'{n}}ski, Miasojedow and
Niemiro}{2012}]{latuszynski-miasojedow-niemiro}
%
\begin{bmisc}[author]
\bauthor{\bsnm{{\L}atuszy{\'{n}}ski},~\bfnm{Krzysztof}\binits{K.}},
\bauthor{\bsnm{Miasojedow},~\bfnm{Blazej}\binits{B.}} \AND
\bauthor{\bsnm{Niemiro},~\bfnm{Wojciech}\binits{W.}}
(\byear{2012}).
\bhowpublished{Nonasymptotic bounds on the estimation error of MCMC algorithms.
\textit{Bernoulli}. To appear.}
\bptok{imsref}%
\end{bmisc}
%
\endbibitem

%b28 #&#
\bibitem[\protect\citeauthoryear{{\L}atuszy{\'n}ski and
Niemiro}{2011}]{latuszynski-niemiro}
%
\begin{barticle}[mr]
\bauthor{\bsnm{{\L}atuszy{\'n}ski},~\bfnm{Krzysztof}\binits{K.}}
\AND
\bauthor{\bsnm{Niemiro},~\bfnm{Wojciech}\binits{W.}}
(\byear{2011}).
\btitle{Rigorous confidence bounds for {MCMC} under a geometric drift
condition}.
\bjournal{J. Complexity}
\bvolume{27}
\bpages{23--38}.
\bid{doi={10.1016/j.jco.2010.07.003}, issn={0885-064X}, mr={2745298}}
\bptok{imsref}%
\end{barticle}
%
\endbibitem

%b29 #&#
\bibitem[\protect\citeauthoryear{Maigret}{1978}]{maigret}
%
\begin{barticle}[mr]
\bauthor{\bsnm{Maigret},~\bfnm{Nelly}\binits{N.}}
(\byear{1978}).
\btitle{Th\'eor\`eme de limite centrale fonctionnel pour une cha\^\i
ne de
{M}arkov r\'ecurrente au sens de {H}arris et positive}.
\bjournal{Ann. Inst. H. Poincar\'e Sect. B (N.S.)}
\bvolume{14}
\bpages{425--440}.
\bid{issn={0020-2347}, mr={0523221}}
\bptok{imsref}%
\end{barticle}
%
\endbibitem

%b30 #&#
\bibitem[\protect\citeauthoryear{Mengersen and
Tweedie}{1996}]{mengersen-tweedie}
%
\begin{barticle}[mr]
\bauthor{\bsnm{Mengersen},~\bfnm{K.~L.}\binits{K.~L.}} \AND
\bauthor{\bsnm{Tweedie},~\bfnm{R.~L.}\binits{R.~L.}}
(\byear{1996}).
\btitle{Rates of convergence of the {H}astings and {M}etropolis algorithms}.
\bjournal{Ann. Statist.}
\bvolume{24}
\bpages{101--121}.
\bid{doi={10.1214/aos/1033066201}, issn={0090-5364}, mr={1389882}}
\bptok{imsref}%
\end{barticle}
%
\endbibitem

%b31 #&#
\bibitem[\protect\citeauthoryear{Metropolis et~al.}{1953}]{metropolis}
%
\begin{barticle}[author]
\bauthor{\bsnm{Metropolis},~\bfnm{N.}\binits{N.}},
\bauthor{\bsnm{Rosenbluth},~\bfnm{A.~W.}\binits{A.~W.}},
\bauthor{\bsnm{Rosenbluth},~\bfnm{M.~N.}\binits{M.~N.}},
\bauthor{\bsnm{Teller},~\bfnm{A.~H.}\binits{A.~H.}} \AND
\bauthor{\bsnm{Teller},~\bfnm{E.}\binits{E.}}
(\byear{1953}).
\btitle{{E}quation of state calculations by fast computing machines}.
\bjournal{J.~Chem. Phys.}
\bvolume{31}
\bpages{1087--1092}.
\bptok{imsref}%
\end{barticle}
%
\endbibitem

%b32 #&#
\bibitem[\protect\citeauthoryear{Meyn and Tweedie}{2009}]{meyn-tweedie}
%
\begin{bbook}[mr]
\bauthor{\bsnm{Meyn},~\bfnm{Sean}\binits{S.}} \AND
\bauthor{\bsnm{Tweedie},~\bfnm{Richard~L.}\binits{R.~L.}}
(\byear{2009}).
\btitle{Markov Chains and Stochastic Stability},
\bedition{2nd} ed.
\bpublisher{Cambridge Univ. Press}, \blocation{Cambridge}.
\bid{mr={2509253}}
\bptok{imsref}%
\end{bbook}
%
\endbibitem

%b33 #&#
\bibitem[\protect\citeauthoryear{Nummelin}{1984}]{nummelin}
%
\begin{bbook}[mr]
\bauthor{\bsnm{Nummelin},~\bfnm{Esa}\binits{E.}}
(\byear{1984}).
\btitle{General Irreducible {M}arkov Chains and Nonnegative Operators}.
\bseries{Cambridge Tracts in Mathematics}
\bvolume{83}.
\bpublisher{Cambridge Univ. Press}, \blocation{Cambridge}.
\bid{doi={10.1017/CBO9780511526237}, mr={0776608}}
\bptok{imsref}%
\end{bbook}
%
\endbibitem

%b34 #&#
\bibitem[\protect\citeauthoryear{Papaspiliopoulos, Roberts and
Sk{\"o}ld}{2007}]{papaspiliopoulos-roberts-skold}
%
\begin{barticle}[mr]
\bauthor{\bsnm{Papaspiliopoulos},~\bfnm{Omiros}\binits{O.}},
\bauthor{\bsnm{Roberts},~\bfnm{Gareth~O.}\binits{G.~O.}} \AND
\bauthor{\bsnm{Sk{\"o}ld},~\bfnm{Martin}\binits{M.}}
(\byear{2007}).
\btitle{A general framework for the parametrization of hierarchical models}.
\bjournal{Statist. Sci.}
\bvolume{22}
\bpages{59--73}.
\bid{doi={10.1214/088342307000000014}, issn={0883-4237}, mr={2408661}}
\bptok{imsref}%
\end{barticle}
%
\endbibitem

%b35 #&#
\bibitem[\protect\citeauthoryear{Papaspiliopoulos and
Roberts}{2008}]{papaspiliopoulos-roberts}
%
\begin{barticle}[mr]
\bauthor{\bsnm{Papaspiliopoulos},~\bfnm{Omiros}\binits{O.}} \AND
\bauthor{\bsnm{Roberts},~\bfnm{Gareth}\binits{G.}}
(\byear{2008}).
\btitle{Stability of the {G}ibbs sampler for {B}ayesian hierarchical models}.
\bjournal{Ann. Statist.}
\bvolume{36}
\bpages{95--117}.
\bid{doi={10.1214/009053607000000749}, issn={0090-5364}, mr={2387965}}
\bptok{imsref}%
\end{barticle}
%
\endbibitem

%b36 #&#
\bibitem[\protect\citeauthoryear{Roberts and
Rosenthal}{1997}]{roberts-rosenthal}
%
\begin{barticle}[mr]
\bauthor{\bsnm{Roberts},~\bfnm{Gareth~O.}\binits{G.~O.}} \AND
\bauthor{\bsnm{Rosenthal},~\bfnm{Jeffrey~S.}\binits{J.~S.}}
(\byear{1997}).
\btitle{Geometric ergodicity and hybrid {M}arkov chains}.
\bjournal{Electron. Commun. Probab.}
\bvolume{2}
\bpages{13--25 (electronic)}.
\bid{doi={10.1214/ECP.v2-981}, issn={1083-589X}, mr={1448322}}
\bptok{imsref}%
\end{barticle}
%
\endbibitem

%b37 #&#
\bibitem[\protect\citeauthoryear{Roberts and
Rosenthal}{2004}]{roberts-rosenthal-survey}
%
\begin{barticle}[mr]
\bauthor{\bsnm{Roberts},~\bfnm{Gareth~O.}\binits{G.~O.}} \AND
\bauthor{\bsnm{Rosenthal},~\bfnm{Jeffrey~S.}\binits{J.~S.}}
(\byear{2004}).
\btitle{General state space {M}arkov chains and {MCMC} algorithms}.
\bjournal{Probab. Surv.}
\bvolume{1}
\bpages{20--71}.
\bid{doi={10.1214/154957804100000024}, issn={1549-5787}, mr={2095565}}
\bptok{imsref}%
\end{barticle}
%
\endbibitem

%b38 #&#
\bibitem[\protect\citeauthoryear{Roberts and Sahu}{1997}]{roberts-sahu}
%
\begin{barticle}[mr]
\bauthor{\bsnm{Roberts},~\bfnm{G.~O.}\binits{G.~O.}} \AND
\bauthor{\bsnm{Sahu},~\bfnm{S.~K.}\binits{S.~K.}}
(\byear{1997}).
\btitle{Updating schemes, correlation structure, blocking and parameterization
for the {G}ibbs sampler}.
\bjournal{J. Roy. Statist. Soc. Ser. B}
\bvolume{59}
\bpages{291--317}.
\bid{doi={10.1111/1467-9868.00070}, issn={0035-9246}, mr={1440584}}
\bptok{imsref}%
\end{barticle}
%
\endbibitem

%b39 #&#
\bibitem[\protect\citeauthoryear{Roberts and Tweedie}{1996}]{roberts-tweedie}
%
\begin{barticle}[mr]
\bauthor{\bsnm{Roberts},~\bfnm{G.~O.}\binits{G.~O.}} \AND
\bauthor{\bsnm{Tweedie},~\bfnm{R.~L.}\binits{R.~L.}}
(\byear{1996}).
\btitle{Geometric convergence and central limit theorems for multidimensional
{H}astings and {M}etropolis algorithms}.
\bjournal{Biometrika}
\bvolume{83}
\bpages{95--110}.
\bid{doi={10.1093/biomet/83.1.95}, issn={0006-3444}, mr={1399158}}
\bptok{imsref}%
\end{barticle}
%
\endbibitem

%b40 #&#
\bibitem[\protect\citeauthoryear{Rockafellar and
Wets}{1998}]{rockafellar-wets}
%
\begin{bbook}[mr]
\bauthor{\bsnm{Rockafellar},~\bfnm{R.~Tyrrell}\binits{R.~T.}} \AND
\bauthor{\bsnm{Wets},~\bfnm{Roger J.~B.}\binits{R.~J.~B.}}
(\byear{1998}).
\btitle{Variational Analysis}.
\bseries{Grundlehren der Mathematischen Wissenschaften [Fundamental Principles
of Mathematical Sciences]}
\bvolume{317}.
\bpublisher{Springer}, \blocation{Berlin}.
\bid{doi={10.1007/978-3-642-02431-3}, mr={1491362}}
\bptok{imsref}%
\end{bbook}
%
\endbibitem

%b41 #&#
\bibitem[\protect\citeauthoryear{Rosenthal}{1995a}]{rosenthal-james-stein}
%
\begin{barticle}[author]
\bauthor{\bsnm{Rosenthal},~\bfnm{Jeffrey~S.}\binits{J.~S.}}
(\byear{1995}a).
\btitle{Analysis of the Gibbs sampler for a model related to James--Stein
estimators}.
\bjournal{Stat. Comput.}
\bvolume{6}
\bpages{269--275}.
\bptok{imsref}%
\end{barticle}
%
\endbibitem

%b42 #&#
\bibitem[\protect\citeauthoryear{Rosenthal}{1995b}]{rosenthal}
%
\begin{barticle}[mr]
\bauthor{\bsnm{Rosenthal},~\bfnm{Jeffrey~S.}\binits{J.~S.}}
(\byear{1995}b).
\btitle{Minorization conditions and convergence rates for {M}arkov chain
{M}onte {C}arlo}.
\bjournal{J. Amer. Statist. Assoc.}
\bvolume{90}
\bpages{558--566}.
\bid{issn={0162-1459}, mr={1340509}}
\bptok{imsref}%
\end{barticle}
%
\endbibitem

%b43 #&#
\bibitem[\protect\citeauthoryear{Roy and Hobert}{2007}]{roy-hobert}
%
\begin{barticle}[mr]
\bauthor{\bsnm{Roy},~\bfnm{Vivekananda}\binits{V.}} \AND
\bauthor{\bsnm{Hobert},~\bfnm{James~P.}\binits{J.~P.}}
(\byear{2007}).
\btitle{Convergence rates and asymptotic standard errors for {M}arkov chain
{M}onte {C}arlo algorithms for {B}ayesian probit regression}.
\bjournal{J. R. Stat. Soc. Ser. B Stat. Methodol.}
\bvolume{69}
\bpages{607--623}.
\bid{doi={10.1111/j.1467-9868.2007.00602.x}, issn={1369-7412}, mr={2370071}}
\bptok{imsref}%
\end{barticle}
%
\endbibitem

%b44 #&#
\bibitem[\protect\citeauthoryear{Stromberg}{1981}]{stromberg}
%
\begin{bbook}[mr]
\bauthor{\bsnm{Stromberg},~\bfnm{Karl~R.}\binits{K.~R.}}
(\byear{1981}).
\btitle{Introduction to Classical Real Analysis}.
\bpublisher{Wadsworth International}, \blocation{Belmont, CA}.
\bid{mr={0604364}}
\bptok{imsref}%
\end{bbook}
%
\endbibitem

%b45 #&#
\bibitem[\protect\citeauthoryear{Tan and Hobert}{2009}]{tan-hobert}
%
\begin{barticle}[mr]
\bauthor{\bsnm{Tan},~\bfnm{Aixin}\binits{A.}} \AND
\bauthor{\bsnm{Hobert},~\bfnm{James~P.}\binits{J.~P.}}
(\byear{2009}).
\btitle{Block {G}ibbs sampling for {B}ayesian random effects models with
improper priors: Convergence and regeneration}.
\bjournal{J. Comput. Graph. Statist.}
\bvolume{18}
\bpages{861--878}.
\bid{doi={10.1198/jcgs.2009.08153}, issn={1061-8600}, mr={2598033}}
\bptok{imsref}%
\end{barticle}
%
\endbibitem

%b46 #&#
\bibitem[\protect\citeauthoryear{Tanner and Wong}{1987}]{tanner-wong}
%
\begin{barticle}[mr]
\bauthor{\bsnm{Tanner},~\bfnm{Martin~A.}\binits{M.~A.}} \AND
\bauthor{\bsnm{Wong},~\bfnm{Wing~Hung}\binits{W.~H.}}
(\byear{1987}).
\btitle{The calculation of posterior distributions by data augmentation}.
\bjournal{J. Amer. Statist. Assoc.}
\bvolume{82}
\bpages{528--550}.
\bid{issn={0162-1459}, mr={0898357}}
\bptok{imsref}%
\end{barticle}
%
\endbibitem

%b47 #&#
\bibitem[\protect\citeauthoryear{Tierney}{1994}]{tierney}
%
\begin{barticle}[mr]
\bauthor{\bsnm{Tierney},~\bfnm{Luke}\binits{L.}}
(\byear{1994}).
\btitle{Markov chains for exploring posterior distributions}.
\bjournal{Ann. Statist.}
\bvolume{22}
\bpages{1701--1762}.
\bid{doi={10.1214/aos/1176325750}, issn={0090-5364}, mr={1329166}}
\bptok{imsref}%
\end{barticle}
%
\endbibitem

\end{thebibliography}
\end{document}

\begin{frontmatter}

\title{Correction\vspace*{6.6pt}\\
Variable transformation to obtain geometric ergodicity
in the random-walk Metropolis algorithm}\vspace*{6pt}
\runtitle{Correction}
\pdftitle{Correction
Variable transformation to obtain geometric ergodicity
in the random-walk Metropolis algorithm}

\textit{Ann. Statist.} \textbf{40} (2013) 3050--3076

\begin{aug}
\author[A]{\fnms{Leif T.} \snm{Johnson}\corref{}\ead[label=e1]{ltjohnson@google.com}}
\and
\author[B]{\fnms{Charles J.} \snm{Geyer}\ead[label=e2]{charlie@stat.umn.edu}}
\runauthor{L. T. Johnson and C. J. Geyer}
\affiliation{Google Inc. and University of Minnesota}
\address[A]{Google Inc.\\
1600 Amphitheatre Parkway \\
Mountain View, California 94043 \\
USA\\
\printead{e1}}
\address[B]{University of Minnesota\\
School of Statistics\\
313 Ford Hall \\
224 Church St. SE \\
Minneapolis, Minnesota 55455\\
USA\\
\printead{e2}}
\end{aug}

% HISTORY:
\received{\smonth{7} \syear{2013}}

% ABSTRACT

% KEYWORDS

\end{frontmatter}

The last sentence of Section~3.1 is incorrect and cannot be fixed
(the rest of Section~3.1 is correct, as is Section~3.2 which uses it).
It is not true that, in general, the random-walk Metropolis algorithm
for Bayesian inference in exponential families with conjugate priors
on natural parameters satisfies the Jarner--Hansen curvature condition.

Poisson regression provides a simple counterexample. Consider the log
likelihood for two independent Poisson random variables $y_1$ and $y_2$
\[
l(\theta)
=
y_1 \theta_1 - e^{\theta_1}
+
y_2 \theta_2 - e^{\theta_2}.
\]
We assume $y_1 > 0$ and $y_2 > 0$ and use flat priors
so the log likelihood is the log unnormalized posterior
and the log unnormalized density of the equilibrium distribution of
the Markov chain.

Consider a sequence $\theta_n$ with components
$\theta_{1, n} \to\infty$, with $\theta_{1, n} > 0$ for all $n$,
and $\theta_{2, n} = - \theta_{1, n}^2$.
It is easily seen that
\[
% \frac{\nabla l(\theta_n)}{\norm{\nabla l(\theta_n)}}
% \to\begin{pmatrix} 1 \\ 0 \end{pmatrix}
% \qquad
% \frac{\theta_n}{\norm{\theta_n}}
% \to\begin{pmatrix} 0 \\ -1 \end{pmatrix}
% \qquad
\frac{\nabla l(\theta_n)}{\norm{\nabla l(\theta_n)}} \cdot
\frac{\theta_n}{\norm{\theta_n}}
\to
0,
\]
where the dot indicates inner product and $\norm{\fatdot}$ denotes the
Euclidean norm. And this shows the Jarner--Hansen curvature condition
fails for this model.

Since this counterexample is a special case of Poisson regression,
we see that geometrically ergodic Markov chain Monte Carlo for
Bayesian inference for Poisson regression with conjugate
priors on natural parameters remains an open research question.


%suskaldyti doi

\printaddresses